\documentclass[a4paper,11pt,reqno]{amsart}

\usepackage{latexsym,dsfont,amssymb,amsmath,amsthm,amscd}

\usepackage[all]{xy}

\chardef\bslash=`\\ 

\numberwithin{equation}{section}

\newtheorem{theorem}{Theorem}[section]
\newtheorem{corollary}[theorem]{Corollary}
\newtheorem{lemma}[theorem]{Lemma}
\newtheorem{proposition}[theorem]{Proposition}
\theoremstyle{remark}
\newtheorem{remark}[theorem]{Remark}

\theoremstyle{definition}

\newcommand\bp{\begin{proof}}
\newcommand\ep{\end{proof}}

\newcommand\unit{{\mathbf 1}}

\newcommand\Dhat{{\hat\Delta}}

\newcommand\CC{{\mathcal C}}

\newcommand\E{{\mathcal E}}
\newcommand\F{{\mathcal F}}
\newcommand\G{{\mathcal G}}

\newcommand\RR{{\mathcal R}}
\newcommand\U{{\mathcal U}}

\newcommand\g{{\mathfrak g}}

\newcommand\Ad{\operatorname{Ad}}
\newcommand\Aut{\operatorname{Aut}}
\newcommand\Inn{\operatorname{Inn}}
\newcommand\Out{\operatorname{Out}}

\newcommand\Hilb{\operatorname{Hilb}_f}
\newcommand\Hom{\operatorname{Hom}}

\newcommand\Ind{\operatorname{Ind}}

\newcommand\Rep{\operatorname{Rep}}
\newcommand\SU{\operatorname{SU}}

\newcommand\Tr{\operatorname{Tr}}

\newcommand{\C}{{\mathbb C}}

\newcommand\Sp{\operatorname{Sp}}
\newcommand\T{{\mathbb T}}
\newcommand\Z{{\mathbb Z}}

\newcommand\eps{\varepsilon}

\newcommand\enu[1]{\smallskip\newline\makebox[5mm][l]{\rm(#1)}}

\begin{document}

\title[Second cohomology]
{On second cohomology of duals of compact groups}


\author[S. Neshveyev]{Sergey Neshveyev}
\address{Department of Mathematics, University of Oslo,
P.O. Box 1053 Blindern, NO-0316 Oslo, Norway
}

\email{sergeyn@math.uio.no}

\author[L. Tuset]{Lars Tuset}
\address{Faculty of Engineering, Oslo University College,
P.O. Box 4 St.~Olavs plass, NO-0130 Oslo, Norway}
\email{Lars.Tuset@iu.hio.no}

\thanks{Supported by the Research Council of Norway.}

\date{November 20, 2010; minor changes May 14, 2011}

\begin{abstract}
We show that for any compact connected group $G$ the second cohomology group defined by unitary invariant $2$-cocycles on $\hat G$ is canonically isomorphic to $H^2(\widehat{Z(G)};\T)$. This implies that the group of autoequivalences of the C$^*$-tensor category $\Rep G$ is isomorphic to $H^2(\widehat{Z(G)};\T)\rtimes\Out(G)$. We also show that a compact connected group $G$ is completely determined by $\Rep G$. More generally, extending a result of Etingof-Gelaki and Izumi-Kosaki we describe all pairs of compact separable monoidally equivalent groups. The proofs rely on the theory of ergodic actions of compact groups developed by Landstad and Wassermann and on its algebraic counterpart developed by Etingof and Gelaki for the classification of triangular semisimple Hopf algebras.

In two appendices we give a self-contained account of amenability of tensor categories, fusion rings and discrete quantum groups, and prove an analogue of Radford's theorem on minimal Hopf subalgebras of quasitriangular Hopf algebras for compact quantum groups.
\end{abstract}

\maketitle

\bigskip

\section*{Introduction}

In an earlier paper \cite{NT4} we showed that any symmetric invariant cocycle on the dual of the $q$-deformation $G_q$ of a simply connected semisimple compact Lie group $G$ is the coboundary of a central element. Our motivation in proving that result was to show that our construction~\cite{NT2} of the Dirac operator on~$G_q$ is canonical up to unitary equivalence. One can, however, interpret this result from a categorical point of view: it shows that any braided autoequivalence of the tensor category~$\Rep G_q$ which is the identity on objects, is naturally isomorphic to the identity functor. Combined with a result of McMullen~\cite{McMul} this implies that the group of braided autoequivalences of~$\Rep G_q$ is isomorphic to the group of automorphisms of the based root datum of $G$. The present paper is motivated by the natural problem of computing the group of all, not necessarily braided, autoequivalences of~$\Rep G_q$. We will deal with  the case $q=1$, for which all the necessary tools have been developed in the study of ergodic actions of compact groups on von Neumann algebras. In a subsequent paper~\cite{NT6} we will extend these results to the quantum case by adapting the technique from~\cite{NT4}.

The notion of a unitary dual $2$-cocycle on a compact group was introduced in the early eighties by Landstad~\cite{La} and Wassermann~\cite{WaPhD,Wa1} in their study of ergodic actions. In Hopf algebra literature it was introduced by Drinfeld~\cite{Dr1}. Some of the results of Landstad and Wassermann were then rediscovered (for finite groups, but non-unitary cocycles) in works of Movshev~\cite{Mo} and Etingof and Gelaki~\cite{EG0,EG1}. Their motivation was however completely different, namely, construction and classification of triangular semisimple Hopf algebras and of fiber functors on their representation categories. The technique was also quite different, with non-unitarity of cocycles imposing additional problems and with a more systematic use of categorical language. The main goal in this paper is to apply the theory developed in~\cite{La,Wa,EG0,EG1} to study invariant dual cocycles on compact connected groups. For finite groups a similar study has been recently undertaken by Guillot and Kassel~\cite{GK}.

In Section~\ref{s1} we collect results on dual $2$-cocycles that will be needed later. Most of them are contained, either explicitly or implicitly, in works of Landstad and Wassermann. To make the paper essentially self-contained we provide complete proofs of most results, relating them to results in Hopf algebra literature and simplifying some of the arguments.

Section~\ref{s2} contains our principal results. We show that the second cohomology group $H^2_G(\hat G;\T)$ defined by unitary invariant dual cocycles is canonically isomorphic to $H^2(\widehat{Z(G)};\T)$ for any compact connected group $G$. There are two facts that make the proof work in the connected group case: a nontrivial one, mentioned earlier, saying that any symmetric invariant cocycle is a coboundary of a central element, and a trivial one, stating that any closed normal abelian subgroup is contained in the center. As an application we show, using a result of McMullen~\cite{McMul}, that the group of autoequivalences of $\Rep G$ is isomorphic to $H^2(\widehat{Z(G)};\T)\rtimes\Out(G)$. Another consequence of the same technique is that a compact connected group is completely determined by the C$^*$-tensor category $\Rep G$ of its finite dimensional unitary representations.

For nonconnected groups the tensor category $\Rep G$ does not always contain full information about~$G$. For finite groups this was shown by Etingof and Gelaki~\cite{EG2}, using the methods they developed in the study of triangular semisimple Hopf algebras, and independently by Izumi and Kosaki~\cite{IK}, using results of Wassermann. In fact, these works describe all pairs of finite groups with equivalent tensor categories of representations. In Section~\ref{s3} we extend this result to compact separable groups. A similar result has been previously obtained by M\"uger, but remained unpublished.

There are two appendices to the paper. In Appendix~\ref{sappa} we give a self-contained presentation of amenability of tensor categories, fusion rings and discrete quantum groups. A simple application we are after is that for a compact group $G$ any unitary fiber functor on $\Rep G$ is dimension preserving, a fact which is used repeatedly in the paper. This is surely known to experts, and amenability is not the only way to prove this, but we use this opportunity to give a unified account of various notions of amenability studied by Longo and Roberts~\cite{LR}, Hiai and Izumi~\cite{HI} and Banica~\cite{Ba1}. In Appendix~\ref{sb} we prove an analogue of a result of Radford \cite{Rad} on minimal Hopf subalgebras of quasitriangular Hopf algebras for compact quantum groups. This is used in Section~\ref{s1} to characterize triangular discrete Hopf $*$-algebras using ideas of Etingof and Gelaki~\cite{EG0,EG1}.

\medskip\noindent
{\bf Acknowledgement.} The first author is grateful to Michael M\"uger, Dmitri Nikshych and Victor Ostrik for stimulating conversations and pointing out some references. Both authors would also like to thank Jurie Conradie and University of Cape Town for hospitality, and Erik B\'edos, Masaki Izumi and Magnus Landstad for comments.

\bigskip

\section{Dual cocycles and ergodic actions} \label{s1}

Let $G$ be a compact group. Denote by $W^*(G)$ the von Neumann algebra of $G$ generated by the operators~$\lambda_g$ of the left regular representation of $G$. This is a coinvolutive Hopf-von Neumann algebra with comultiplication $\Dhat(\lambda_g)=\lambda_g\otimes\lambda_g$ and antipode $\hat S(\lambda_g)=\lambda_{g^{-1}}$. Denote also by $\hat\eps$ the counit of~$W^*(G)$. Denote by $\U(G)$ the algebra of closed densely defined operators affiliated with~$W^*(G)$. In other words, if we choose representatives $\pi_j\colon G\to B(H_j)$ of the isomorphism classes of irreducible unitary representations of $G$ and identify~$W^*(G)$ with the $\ell^\infty$-sum of the algebras~$B(H_j)$, then $\U(G)=\prod_jB(H_j)$. Denote by $\Rep G$ the C$^*$-tensor category of finite dimensional unitary representations of $G$, and by $\Hilb$ the C$^*$-tensor category of finite dimensional Hilbert spaces.

\smallskip

An invertible element $\F\in \U(G\times G)$ is called a $2$-cocycle on $\hat G$, or a dual $2$-cocycle on $G$, if
\begin{equation} \label{ecocycle0}
(\F\otimes1)(\Dhat\otimes\iota)(\F)=(1\otimes\F)(\iota\otimes\Dhat)(\F).
\end{equation}
Two cocycles $\F$ and $\G$ are called cohomologous if there exists an invertible element $u\in\U(G)$ such that $\G=(u\otimes u)\F\Dhat(u)^{-1}$; a cocycle of the from $(u\otimes u)\Dhat(u)^{-1}$ is called a coboundary. We denote by $Z^2(\hat G;\C^*)$ the set of $2$-cocycles, and by $H^2(\hat G;\C^*)$ the set of cohomology classes of $2$-cocycles. If~$G$ is abelian, so that $\hat G$ is a discrete abelian group, then $H^2(\hat G;\C^*)$ is the usual second cohomology group.

\smallskip

We will be mainly interested in unitary $2$-cocycles $\F$, so that $\F$ is a unitary element in the von Neumann algebra $W^*(G)\bar\otimes W^*(G)\subset\U(G\times G)$. Correspondingly, two such cocycles $\F$ and $\G$ are called cohomologous if there exists a unitary element $u\in W^*(G)$ such that $\G=(u\otimes u)\F\Dhat(u)^*$, and the set of cohomology classes of unitary $2$-cocycles is denoted by $H^2(\hat G;\T)$.

\smallskip

Note that by applying $\iota\otimes\hat\eps\otimes\iota$ to the cocycle identity~\eqref{ecocycle0}, we  see that
$(\iota\otimes\hat\eps)(\F)\otimes1=1\otimes(\hat\eps\otimes\iota)(\F)$, which implies that
$$
(\iota\otimes\hat\eps)(\F)=(\hat\eps\otimes\iota)(\F)=(\hat\eps\otimes\hat\eps)(\F)1.
$$
Thus replacing $\F$ by the cohomologous cocycle $(\hat\eps\otimes\hat\eps)(\F)^{-1}\F$, we may assume that $(\iota\otimes\hat\eps)(\F)=(\hat\eps\otimes\iota)(\F)=1$. Such cocycles are called counital.

\smallskip

Assume the group $G$ acts on a von Neumann algebra $M$. The action is called ergodic if the fixed point algebra $M^G$ is trivial. If the action is ergodic, it is said to be of full multiplicity if
$$
\dim\Hom_G(H_U,M)=\dim H_U
$$
for every irreducible representation $\pi_U\colon G\to B(H_U)$ of $G$. The equality then clearly holds for all finite dimensional representations. Note that the inequality $\le$ holds for any ergodic action~\cite{HKLS}. The linear span $M_U$ of the spaces $\gamma(H_U)$ for $\gamma\in \Hom_G(H_U,M)$ is called the spectral subspace of $M$ corresponding to $U$.

\smallskip

The main part of the following theorem, the correspondence between (i) and (iii), was proved by Landstad~\cite{La} and Wassermann~\cite{WaPhD,Wa1}.

\begin{theorem} \label{tergact}
For any compact group $G$ there are one-to-one correspondences between
\enu{i} isomorphism classes of full multiplicity ergodic actions of $G$ on von Neumann algebras;
\enu{ii} unitary isomorphism classes of unitary fiber functors $\Rep G\to\Hilb$;
\enu{iii} elements of $H^2(\hat G;\T)$.
\end{theorem}

\bp For the reader's convenience we give a proof, referring to \cite{Wa1} for a more detailed although somewhat different argument. Yet another proof of a more general result valid for compact quantum groups can be found in~\cite{BDRV}. Even more general result is proved in~\cite{PR}, where ergodic actions of compact quantum groups are described in terms of weak tensor functors.

\smallskip

Assume $\alpha\colon G\to \Aut(M)$ is an ergodic action. Define a functor $F\colon \Rep G\to\Hilb$ by letting
$$
F(U)=(M\otimes H_U)^G\cong\Hom_G(H_U^*,M)
$$
for any finite dimensional unitary representation $U\colon G\to B(H_U)$ of $G$. The action of $F$ on morphisms is defined in the obvious way. The scalar product on $F(U)$ is defined as the restriction of that on~$M\otimes H_U$ defined by
$$
(a\otimes\xi,b\otimes\zeta)=\varphi(b^*a)(\xi,\zeta),
$$
where $\varphi$ is the unique (thanks to ergodicity) $G$-invariant normal state on $M$. Note that if $X=\sum_ia_i\otimes\xi_i$ and $Y=\sum_jb_j\otimes\zeta_j$ are in $(M\otimes H_U)^G$ then $(X,Y)1=\sum_{i,j}(\xi_i,\zeta_j)b_j^*a_i$, since the latter expression is $G$-invariant. By slightly abusing notation this can be written as $(X,Y)1=Y^*X$. The opposite multiplication map $M\otimes_{alg}M\to M$, $a\otimes b\mapsto ba$, induces a map
$$
F_2(U,V)\colon F(U)\otimes F(V)\to F(U\otimes V).
$$
Using that the scalar products are defined by $(X,Y)1=Y^*X$, it is easy to check that this is an isometry, hence it is unitary by the full multiplicity assumption. Thus $(F,F_2)$ is a unitary fiber functor.

\smallskip

Assume now that we are given a unitary fiber functor $F\colon \Rep G\to\Hilb$. By Corollary~\ref{cdimpres} below, it preserves the dimensions of objects. Since $\Rep G$ is semisimple, it follows that disregarding the tensor structure, the functor $F$ is isomorphic to the forgetful one. Choose natural unitary isomorphisms $\eta_U\colon H_U\to F(U)$. Then we get natural unitaries $\eta^{-1}_{U\otimes V}F_2(U,V)(\eta_U\otimes\eta_V)$ on $H_U\otimes H_V$. By naturality they are defined by the action of a unitary $\F^*\in W^*(G)\bar\otimes W^*(G)$. It is easy to see that $\F$ is a unitary $2$-cocycle, that a different choice of isomorphisms $\eta_V$ would give a cohomologous cocycle, and that up to a natural unitary  isomorphism the fiber functor $F$ is completely determined by $\F$; see e.g.~\cite[Proposition 1.1]{NT3} for details.

\smallskip

Finally, assume we are given a unitary $2$-cocycle $\F\in Z^2(\hat G;\T)$. Let $\C[G]\subset C(G)$ be the algebra of regular functions on $G$. It can be considered as a subalgebra of the predual $W^*(G)_*$ of $W^*(G)$: if $a\in \C[G]$ then $a(\lambda_g)=a(g)$. Using this define a new algebra~$\C_\F[G]$ with the same underlying space~$\C[G]$ and new product
$$
ab=(a\otimes b)(\F\Dhat(\cdot)).
$$
In other words, $ab=(a_{(0)}\otimes b_{(0)})(\F)a_{(1)}\cdot b_{(1)}$, where $\cdot$ denotes the usual pointwise product in $\C[G]\subset C(G)$. Note that the unit of $\C_\F[G]$ is the element ${(\hat\eps\otimes\hat\eps)(\F)}^{-1}1\in\C[G]$. We will see soon that there is a canonical $*$-structure on $\C_\F[G]$. Define a representation $\pi_\F$ of $\C_\F[G]$ on $L^2(G)$ as follows. Let $W\in B(L^2(G)\otimes L^2(G))$ be the fundamental multiplicative unitary, so $(W\xi)(g,h)=\xi(g,g^{-1}h)$. It has the properties
$$
W\in L^\infty(G)\bar\otimes W^*(G), \ \ \Dhat(x)=W(x\otimes1)W^*\ \ \hbox{for} \ \ x\in W^*(G),\ \ (\iota\otimes\Dhat)(W)=W_{12}W_{13}.
$$
Using the last two properties we get
\begin{equation} \label{ecocycle}
(\F W)_{12}(\F W)_{13}=\F_{23}(\iota\otimes\Dhat)(\F W)\ \ \hbox{in}\ \ B(L^2(G))\bar\otimes W^*(G)\bar\otimes W^*(G),
\end{equation}
since $(\F W)_{12}(\F W)_{13}=\F_{12}(\Dhat\otimes\iota)(\F)W_{12}W_{13}$ and $\F_{23}(\iota\otimes\Dhat)(\F W)
=\F_{23}(\iota\otimes\Dhat)(\F)W_{12}W_{13}$. This implies that the formula $\pi_\F(a)=(\iota\otimes a)(\F W)$ defines a unital representation of $\C_\F[G]$ on~$L^2(G)$. Denote by $A$ the norm closure of the image of $\pi_\F$. We want to show that this is a C$^*$-algebra. Denoting by $\rho_g$ the operators of the right regular representation of $G$, we have
\begin{equation}\label{eaction}
(\rho_g\otimes1)W(\rho_g^*\otimes1)=W(1\otimes\lambda_g).
\end{equation}
Since the operators $\rho_g\otimes1$ commute with $\F\in W^*(G)\bar\otimes W^*(G)$, we have a similar identity for $\F W$ instead of $W$. It follows that we can define an action $\alpha$ of $G$ on $A$ by letting $\alpha_g=\Ad\rho_g|_A$. Furthermore, for $a\in\C_\F[G]=\C[G]$ we have $\alpha_g(\pi_\F(a))=\pi_\F(a(\cdot\, g))$. Since every orbit of the action of $G$ on $\pi_\F(\C_\F[G])$ lies in a finite dimensional space, the action of $G$ on $A$ is continuous. As the action of $G$ on itself by right translations is ergodic, we have $\int_G\alpha_g(a)dg\in\C1$ for every $a\in\pi_\F(\C_\F[G])$, hence for every $a\in A$. Therefore the action of~$G$ on~$A$ is ergodic. Using~\eqref{eaction} once again we see that for every $\omega\in\C[G]\subset W^*(G)_*$ the operator
$$
(\iota\otimes\omega)(\F W(\iota\otimes\hat S)(\F W))
$$
is a $G$-invariant element of $A$, hence a scalar. Therefore $\F W(\iota\otimes\hat S)(\F W)\in1\otimes \U(G)$, so there exists an element $u\in \U(G)$ such that
\begin{equation} \label{edual}
(\iota\otimes\hat S)((\F W)^*)=(1\otimes u)\F W.
\end{equation}
This implies that the algebra $\pi_\F(\C_\F[G])$ is self-adjoint, so $A$ is indeed a C$^*$-algebra. We will denote it by $C^*(\hat G;\F)$. Taking the weak operator closure $M$ of $A$, denoted also by $W^*(\hat G;\F)$, we get a von Neumann algebra with an ergodic action of $G$, which we continue to denote by $\alpha$. It is easy to check that if we take a cohomologous cocycle then we get an isomorphic ergodic action.

We still have to show that the action $\alpha$ on $M=W^*(\hat G;\F)$ has full multiplicity. For every finite dimensional unitary representation $U\colon G\to B(H_U)$ put
$$
\Theta_U=(\iota\otimes\pi_U)(\F W)\in M\otimes B(H_U),
$$
where $\pi_U$ is the extension of $U$ to $W^*(G)$. Then $(\alpha_g\otimes\iota)(\Theta_U)=\Theta_U(1\otimes\pi_U(g))$ by~\eqref{eaction}. It follows that the map
$$
\eta_U\colon H_U\to M\otimes H_U,\ \ \xi\mapsto \Theta_U^*(1\otimes\xi),
$$
is an isometric embedding of $H_U$ into $(M\otimes H_U)^G\cong \Hom_G(H_U^*,M)$. Since the dimension of the latter space is not bigger than $\dim H_U$ by ergodicity, it follows that $\eta_U\colon H_U\to(M\otimes H_U)^G$ is an isomorphism. In particular, the ergodic action $\alpha$ of~$G$ on~$M$ has full multiplicity. It follows that since the spectral subspaces of the action of~$G$ on~$M$ are contained in $\pi_\F(\C_\F[G])$, the representation~$\pi_\F$ is faithful on~$\C_\F[G]$. Hence $\C_\F[G]$ inherits a $*$-structure from $A$; explicitly, $a^*(x)=\overline{a(\hat S(ux)^*)}$ for $x\in W^*(G)$, where $u$ is given by~\eqref{edual}. In other words, $a^*=a^\dag_{(0)}(u)a^\dag_{(1)}$, where $\dag$ denotes the usual $*$-operation on $\C[G]\subset C(G)$.

\smallskip

To complete the proof of the theorem we have to show that any full multiplicity ergodic action is isomorphic to the action on $W^*(\hat G;\F)$ for some cocycle $\F$, and that the cocycle defined by the action on $W^*(\hat G;\F)$ is cohomologous to $\F$. We start with the second problem.

\smallskip

Let $M=W^*(\hat G;\F)$. As above, consider the isomorphisms $\eta_U(\xi)=\Theta_U^*(1\otimes\xi)$ with $\Theta_U=(\iota\otimes\pi_U)(\F W)$. Then identity~\eqref{ecocycle} can be written as
\begin{equation} \label{ecocycle2}
1\otimes(\pi_{U}\otimes\pi_{V})(\F^*)=\Theta_{U\otimes V}(\Theta_U)^*_{13}(\Theta_V)^*_{12}\ \ \hbox{in}\ \ M\otimes B(H_U)\otimes B(H_V).
\end{equation}
But this exactly means that $\F$ is the cocycle defined by the ergodic action on $M$ and the isomorphisms~$\eta_U$.

\smallskip

Assume now that we have a full multiplicity ergodic action~$\alpha$ on a von Neumann algebra $M$. Choose natural unitary isomorphisms $\eta_U\colon H_U\to (M\otimes H_U)^G$. Let $\Theta_U\in M\otimes B(H_U)$ be such that $\eta_U(\xi)=\Theta_U^*(1\otimes\xi)$. Since $\eta_U$ is unitary, we have $(\varphi\otimes\iota)(\Theta_U\Theta_U^*)=1$, where $\varphi$ is the unique normal $G$-invariant state on~$M$. On the other hand, since the image of $\eta_U$ is $G$-invariant, we have $(\alpha_g\otimes\iota)(\Theta_U)=\Theta_U(1\otimes U(g))$, whence $\Theta_U\Theta_U^*\in 1\otimes B(H_U)$ by ergodicity. Hence $\Theta_U\Theta_U^*=1$, and since the algebra~$M$ is finite by~\cite[Theorem~4.1]{HKLS}, we conclude that also $\Theta_U^*\Theta_U=1$. By naturality there exists a unitary $\Theta\in M\bar\otimes W^*(G)$ such that $(\iota\otimes\pi_U)(\Theta)=\Theta_U$ for all $U$, and then $(\alpha_g\otimes\iota)(\Theta)=\Theta(1\otimes\lambda_g)$.

Let $\F$ be the cocycle defined by the action $\alpha$ and the isomorphisms $\eta_U$. Then $\F$ satisfies identity~\eqref{ecocycle2}, from which we conclude that
\begin{equation} \label{ecocycle3}
\Theta_{12}\Theta_{13}=\F_{23}(\iota\otimes\Dhat)(\Theta).
\end{equation}
Let $u\in \U(G)$ be the element defined by~\eqref{edual}. We claim that then
\begin{equation} \label{edual2}
(\iota\otimes\hat S)(\Theta^*)=(1\otimes u)\Theta.
\end{equation}
To show this, consider the unitary $\tilde\Theta=\Theta_{13}W_{23}\in M\bar\otimes B(L^2(G))\bar\otimes W^*(G)$. Using~\eqref{ecocycle3} and that $\Dhat(x)W=W(x\otimes1)$, we get $\Theta_{12}\tilde\Theta\Theta_{12}^*=(\F W)_{23}$, whence
$$
(\iota\otimes\iota\otimes\hat S)(\tilde\Theta^*)=\Theta_{12}^*(\iota\otimes\iota\otimes\hat S)((\F W)_{23}^*)\Theta_{12}
=(1\otimes1\otimes u)\tilde\Theta.
$$
On the other hand, since $(\iota\otimes\hat S)(W^*)=W$, we have
$$
(\iota\otimes\iota\otimes\hat S)(\tilde\Theta^*)=(\iota\otimes\iota\otimes\hat S)(\Theta^*_{13})W_{23}
=(\iota\otimes\iota\otimes\hat S)(\Theta^*_{13})\Theta_{13}^*\tilde\Theta,
$$
and~\eqref{edual2} is proved.

Identities~\eqref{ecocycle3} and~\eqref{edual2} imply that by letting $\pi(a)=(\iota\otimes a)(\Theta)$ for $a\in\C_\F[G]$ we get a $G$-equivariant unital $*$-homomorphism $\pi\colon \C_\F[G]\to M$. Then $\varphi\circ \pi$ is a $G$-invariant functional on $\C_\F[G]$ with value~$1$ at the unit. But by ergodicity there exists only one such functional, the restriction of the unique normal $G$-invariant state on $W^*(\hat G;\F)$ to $ \C_\F[G]$. Hence $\pi$ extends to a normal $G$-equivariant embedding $W^*(\hat G;\F)\to M$. Since the actions of $G$ on both algebras are of full multiplicity, we conclude that this is an isomorphism.
\ep

\begin{remark} \label{rnormalized}
The element $u\in\U(G)$ defined by \eqref{edual} is in fact unitary. Indeed, consider the element $U=(\iota\otimes\hat S)(\F W)$. By ergodicity the element $U^*U$ has the form $1\otimes w$ for some $w\in\U(G)$. Using that the unique $G$-invariant state $\varphi$ on $W^*(\hat G;\F)$ is tracial by~\cite[Theorem~4.1]{HKLS}, we get
$$
w=(\varphi\otimes\iota)(U^*U)=(\varphi\otimes\hat S)(\F W(\F W)^*)=1,
$$
so $U$ is unitary and hence $u$ is unitary as well. Furthermore, by \cite[Lemma 10]{Wa1} the element $u$ satisfies the identity
$$
\Dhat(e_1)(1\otimes u)=\overline{(\hat\eps\otimes\hat\eps)(\F)}\Dhat(e_1)\F^*,
$$
where $e_1\in W^*(G)$ is the central projection corresponding to the trivial representation, so that $ae_1=\hat\eps(a)e_1$. It follows that $u=\overline{(\hat\eps\otimes\hat\eps)(\F)}m(\hat S\otimes\iota)(\F^*)$, where $m\colon W^*(G)\bar\otimes W^*(G)\to\U(G)$ is the multiplication map. Thus $u$ is an element familiar from Hopf algebra theory: if $\Dhat_\F=\F\Dhat(\cdot)\F^*$ is the twisted coproduct, then the new antipode is given by $\hat S_\F=u^*\hat S(\cdot)u$, see e.g.~\cite[Proposition~4.2.13]{CP}. Note that one can easily see from~ \eqref{edual} that $u$ is $\hat S$-invariant,
which implies that $\hat S_\F^2=\iota$. Furthermore, if we replace $\F$ by a cohomologous cocycle $(v\otimes v)\F\Dhat(v)^*$, then $u$ changes to $\hat S(v^*)uv^*$, and from this it is not difficult to see that any cocycle is cohomologous to a counital cocycle such that the corresponding element $u$ equals $1$, see \cite[Theorem~5]{Wa1}, so in this case $\hat S_\F=\hat S$.  Such cocycles are called normalized.

\end{remark}

An important consequence of the above theorem is that if a unitary cocycle $\F$ is symmetric, that is, $\F_{21}=\F$, then it is a coboundary. Indeed, in this case the algebra $W^*(\hat G;\F)$ is abelian. But there exists only one full multiplicity ergodic action of $G$ on an abelian von Neumann algebra, namely, the action of $G$ on $L^\infty(G)$ by right translations, and this action corresponds to the trivial cocycle.

\smallskip

Assume $H\subset G$ is a closed subgroup. Then $W^*(H)$ embeds into $W^*(G)$, and therefore any unitary $2$-cocycle~$\E$ on $\hat H$ can be considered as a cocycle $\F$ on $\hat G$. We will say that $\F$ is induced from $H$ and write $\F=\Ind^G_H\E$. From the point of view of ergodic actions this seemingly trivial construction corresponds to induction, defined as follows. Let $\beta\colon H\to\Aut(N)$ be an ergodic action. Consider the action $\gamma$ of $H$ on $L^\infty(G)\bar\otimes N$ defined by $\gamma_h(f\otimes a)=f(\cdot\,h)\otimes\beta_h(a)$. Let $M$ be the fixed point algebra $(L^\infty(G)\bar\otimes N)^H$. The action of $G$ on $L^\infty(G)$ by left translations gives an ergodic action $\alpha$ of~$G$ on $M$. We write $(M,\alpha)=\Ind^G_H(N,\beta)$. Note that $\alpha$ has full multiplicity if and only if $\beta$ has full multiplicity, which can be easily seen e.g.~using Frobenius reciprocity.

\begin{lemma} \label{lind}
Assume $G$ is a compact group and $H\subset G$ is a closed subgroup.
For $\F\in Z^2(\hat G;\T)$ and the corresponding action $\alpha$ of $G$ on $M=W^*(\hat G;\F)$ the following conditions are equivalent:
\enu{i} there exists a $G$-equivariant embedding of $L^\infty(G/H)$ into the center of $M$;
\enu{ii} $\alpha$ is isomorphic to an action induced from $H$;
\enu{iii} $\F$ is cohomologous to a cocycle induced from $H$.
\end{lemma}

\bp  (i)$\Leftrightarrow$(ii). This is proved in~\cite[Theorem~7]{Wa}. We shall briefly present the argument, since it will be partially used later. The implication (ii)$\Rightarrow$(i) is obvious. To prove the converse, consider the weakly operator dense C$^*$-subalgebra $A\subset M$ of $G$-continuous elements, that is, elements $a\in M$ such that the map $G\ni g\mapsto \alpha_g(a)\in M$ is norm-continuous. By assumption, $A$ is a $C(G/H)$-algebra. Let~$A_e$ be its fiber at $H\in G/H$, and $\pi_e\colon A\to A_e$ be the quotient map. The action of $H$ on $A$ defines an action of $H$ on $A_e$. We have a $G$-equivariant embedding $A\hookrightarrow C(G)\otimes A_e= C(G;A_e)$, $a\mapsto f_a$, $f_a(g)=\pi_e(\alpha_g^{-1}(a))$. It can be checked that it is an isomorphism of $A$ onto $(C(G)\otimes A_e)^H$. In particular, the action of $H$ on $A_e$ is ergodic. Therefore $\alpha$ is induced from the action of $H$ on the von Neumann algebra generated by~$A_e$ in the GNS-representation defined by the unique $H$-invariant state on $A_e$.

\smallskip

(ii)$\Leftrightarrow$(iii). It suffices to show that if $\E\in Z^2(\hat H;\T)$, $\beta$ is the corresponding action of $H$ on $N=W^*(\hat H;\E)$, and $(L,\gamma)=\Ind^G_H(N,\beta)$, then the cocycle defined by the action $\gamma$ is cohomologous to $\Ind^G_H\E$. By the proof of Theorem~\ref{tergact} this is the case if and only if there exists a unitary $\Theta\in L\bar\otimes W^*(G)$ such that
$$
(\gamma_g\otimes\iota)(\Theta)=\Theta(1\otimes\lambda_G(g))\ \ \hbox{and}\ \ \Theta_{12}\Theta_{13}=\E_{23}(\iota\otimes\Dhat)(\Theta),
$$
where $\lambda_G$ denotes the left regular representation of $G$. Let $\Omega\in N\bar\otimes W^*(H)$ be such a unitary for the action $\beta$ of $H$ on $N$, so
$$
(\beta_h\otimes\iota)(\Omega)=\Omega(1\otimes\lambda_H(h))\ \ \hbox{and}\ \ \Omega_{12}\Omega_{13}=\E_{23}(\iota\otimes\Dhat)(\Omega),
$$
where $\lambda_H$ denotes the left regular representation of $H$. Define
$$
\Theta\in L^\infty(G)\bar\otimes N\bar\otimes W^*(G)=L^\infty(G;N\bar\otimes W^*(G))\ \ \hbox{by}\ \ \Theta(g)=\Omega(1\otimes\lambda_G(g^{-1})).
$$
It is straightforward to check that $\Theta$ lies in $(L^\infty(G)\bar\otimes N)^H\bar\otimes W^*(G)=L\bar\otimes W^*(G)$ and has the right properties.
\ep

For any ergodic action of $G$ on $M$, the action on the center $Z(M)$ is ergodic, hence $Z(M)\cong L^\infty(G/H)$ for a closed subgroup $H\subset G$. By the lemma above we conclude that a cocycle $\F$ cannot be cohomologous to a cocycle induced from a proper subgroup if and only if $W^*(\hat G;\F)$ is a factor. Such cocycles are called nondegenerate. As was observed already by Albeverio and H{\o}egh-Krohn in the preliminary version of~\cite{AHK}, the study of ergodic actions is reduced to that of ergodic actions on factors. We therefore have the following result; for finite groups and $\C^*$-valued cocycles it was proved by Movshev~\cite{Mo}, Etingof and Gelaki~\cite{EG1} (see also \cite{CE}) and Davydov~\cite{Da}.

\begin{theorem} \label{tred}
Let $G$ be a compact group. Then for any full multiplicity ergodic action $\alpha$ of $G$ on a von Neumann algebra $M$ there exists a closed subgroup $H\subset G$ and a full multiplicity ergodic action $\beta$ of $H$ on a factor $N$ such that $(M,\alpha)=\Ind^G_H(N,\beta)$. If $(\tilde H,\tilde N,\tilde\beta)$ is another such triple then there exists an element $g\in G$ such that $gHg^{-1}=\tilde H$ and the action $\beta$ on $N$ is isomorphic to the action $h\mapsto\tilde\beta_{ghg^{-1}}$ of $H$ on $\tilde N$.

Equivalently, for any cocycle $\F\in Z^2(\hat G;\T)$ there exists a closed subgroup $H$ of $G$ and a nondegenerate cocycle $\E\in Z^2(\hat H;\T)$ such that $\F$ is cohomologous
to $\Ind^G_H\E$. If $(\tilde H,\tilde\E)$ is another such pair then there exists an element $g\in G$ such that $gHg^{-1}=\tilde H$ and the cocycles $\E$ and $(\lambda_g\otimes\lambda_g)^{-1}\tilde\E(\lambda_g\otimes\lambda_g)$ are cohomologous as cocycles on $\hat H$.
\end{theorem}

\bp As we have already remarked above, the center $Z(M)$ can be identified with $L^\infty(G/H)$. Then by Lemma~\ref{lind} the action $\alpha$ is induced from an action $\beta$ of $H$ on a von Neumann algebra $N$, and $N$ is necessarily a factor, since otherwise the center of $M$, which contains $(L^\infty(G)\bar\otimes Z(N))^H\cong L^\infty(G/H)\bar\otimes Z(N)$, would be strictly larger than $L^\infty(G/H)$.

To show uniqueness, as in the proof of Lemma~\ref{lind} consider the C$^*$-subalgebra $A\subset M$ of $G$-continuous elements. Then if $(M,\alpha)=\Ind^G_H(N,\beta)$ with $N$ a factor, the center of $A$ is canonically $G$-equivariantly isomorphic to $C(G/H)$. Hence if $(M,\alpha)=\Ind^G_{\tilde H}(\tilde N,\tilde \beta)$ with $\tilde N$ a factor for another subgroup $\tilde H$, then the spaces $G/H$ and $G/\tilde H$ are $G$-equivariantly homeomorphic, so $H$ and $\tilde H$ are inner conjugate in $G$. It remains to show that the pair $(N,\beta)$ is uniquely defined for a fixed identification of the center of $M$ with $L^\infty(G/H)$. Let $B$ be the C$^*$-subalgebra $B\subset N$ of $H$-continuous elements. Then by definition of induction, $B$ is $H$-equivariantly isomorphic to the fiber $A_e$ of the $C(G/H)$-algebra $A$ at the point $H\in G/H$.
\ep

The set $H^2(\hat G;\T)$ can therefore be described as follows. Denote by $S(G)$ the set of inner conjugacy classes of closed subgroups of $G$. For every $s\in S(G)$ choose a representative $G_s$. Then
\begin{equation} \label{eh2}
H^2(\hat G;\T)=\bigsqcup_{s\in S(G)}H^2(\hat G_s;\T)^\times/N_G(G_s),
\end{equation}
where $H^2(\hat G_s;\T)^\times\subset H^2(\hat G_s;\T)$ is the subset of elements represented by nondegenerate cocycles, and $N_G(G_s)$ is the normalizer of $G_s$ in $G$. The action of $N_G(G_s)$ on $Z^2(\hat G_s;\T)$ is  defined by $\E\mapsto(\lambda_g\otimes\lambda_g)\E(\lambda_g\otimes\lambda_g)^{-1}$.

For finite groups this can be formulated in a better way~\cite{Mo}. Assume a finite group $G$ acts ergodically with full multiplicity on  a factor $M$. Then $M=B(H)$ for a finite dimensional Hilbert space $H$, and for dimension reasons $\dim H=|G|^{1/2}$. Then the action is implemented by a uniquely defined projective unitary representation $\rho\colon G\to PU(H)$. Lift $\rho$ to a map $\tilde\rho\colon G\to U(H)$ and consider the cocycle $c\in Z^2(G;\T)$ defined by $c(g,h)=\tilde\rho(g)\tilde\rho(h)\tilde\rho(gh)^*$; note that the cohomology class of~$c$ is uniquely defined by our action. Consider the twisted group C$^*$-algebra $C^*(G;c)$ with unitary generators $u_g$, so that $u_gu_h=c(g,h)u_{gh}$. Define a $*$-homomorphism $\pi\colon C^*(G;c)\to B(H)$ by $\pi(u_g)=\tilde\rho(g)$. Then $\pi(C^*(G;c))'=B(H)^G=\C1$, so $\pi\colon C^*(G;c)\to B(H)$ is onto, hence it is an isomorphism for dimension reasons. It follows that the action of $G$ on $B(H)$ we started with is isomorphic to the action of $G$ on $C^*(G;c)$ defined by $g\mapsto\Ad u_g$, so the cocycle~$c$ completely determines the action. Since $C^*(G;c)$ is a factor, the cocycle $c$ is nondegenerate in the usual sense: for every $g\in G$, $g\ne e$, the character of the centralizer $C(g)$ of $g$ in $G$ defined by $C(g)\ni h\mapsto c(g,h)c(h,g)^{-1}$ is nontrivial. Denote by $H^2(G;\T)^\times\subset H^2(G;\T)$ the subset of elements represented by nondegenerate cocycles. We have therefore constructed an injective map $H^2(\hat G;\T)^\times\to H^2(G;\T)^\times$. In fact, it is bijective: if $c\in Z^2(G;\T)$ is a nondegenerate cocycle then $C^*(G;c)$ is a factor, so the action $g\mapsto\Ad u_g$ on it is ergodic, and for dimension reasons it has full multiplicity. Thus for any finite group $G$ equality~\eqref{eh2} gives
\begin{equation} \label{eh2a}
H^2(\hat G;\T)\cong\bigsqcup_{s\in S(G)}H^2(G_s;\T)^\times/N_G(G_s).
\end{equation}

It is clear that even in the finite group case the computation of $H^2(\hat G;\T)$ using~\eqref{eh2} or~\eqref{eh2a} for any large class of groups is a daunting task, as it in particular involves classification up to inner conjugacy of all subgroups of $G$ supporting a nondegenerate cocycle. A finite group for which there is a nondegenerate $2$-cocycle is called a central type factor group (ctfg). We mention a few facts about ctfgs, see e.g.~\cite{Sh} and references therein for more information. One immediate property of ctfgs is that their orders are squares. An abelian group is a ctfg if and only if it has the form~$H\times H$. Any ctfg is solvable; this is a highly nontrivial result established using the classification of finite simple groups. Any finite solvable group can be embedded into a ctfg. The lowest order for which there exists a non-abelian ctfg is $16$, and $\Z/2\Z\times D_8$ is probably the simplest example of such a group, see~\cite{IK2},  pages 146--147. Another concrete example of a non-abelian ctfg is $\Z/2\Z\times(\Z/3\Z\rtimes S_3)$, where~$S_3$ acts on $\Z/3\Z$ by the sign of permutation, see~\cite[Example~4.3]{EG}.

Turning to infinite compact groups, the question when such a group has a nondegenerate dual cocycle, or equivalently, when $G$ has an ergodic action of full multiplicity on an injective II$_1$-factor, is wide open. In view of infinitesimal considerations in~\cite{Wa2} and~\cite{Mo}, a concrete question one might ask is whether a connected compact group with a nondegenerate dual cocycle is necessarily abelian. A related question, posed in~\cite{HKLS}, is whether a compact simple group cannot have an ergodic action on a II$_1$-factor. The only results in this direction are the computations of Wassermann~\cite{Wa2} showing that for the groups $SU(2)\times SU(2)$ and $SU(3)$ and their quotients, any dual cocycle is cohomologous to one induced from a maximal torus. This in particular implies that these groups do not have non-abelian subgroups with nondegenerate cocycles. Since any finite group embeds into $SU(n)$ for sufficiently large $n$, as the rank grows the situation clearly becomes more complicated.

\smallskip

It would also be interesting to give a description of $H^2(\hat G;\C^*)$. For finite groups the answer is similar to~\eqref{eh2a}: $H^2(\hat G;\C^*)\cong\bigsqcup_{s\in S(G)}H^2(G_s;\C^*)^\times/N_G(G_s)$, see~\cite{Mo,EG1,Da}. Since a $\C^*$-valued $2$-cocycle on a finite group $H$ is cohomologous to a cocycle with values in the group of roots of unity of order $|H|$, it follows that the canonical map $H^2(\hat G;\T)\to H^2(\hat G;\C^*)$ is a bijection. This is no longer true already for infinite abelian groups. Furthermore, the structure of non-unitary cocycles is much more complicated. For example, as we mentioned above, any unitary dual cocycle on $SU(3)$ is cohomologous to one  induced from a maximal torus, while the computations of Ohn~\cite{Ohn} show that there are non-unitary cocycles that cannot be obtained this way. This is apparently related to the fact that the complexification~$G_\C$ of a compact Lie group $G$ has many more Poisson-Lie structures than the group~$G$ itself.

\smallskip

For a compact group $G$ we called a cocycle $\F\in Z^2(\hat G;\T)$ nondegenerate if $W^*(\hat G;\F)$ is a factor, or equivalently, if $\F$ is not cohomologous to a cocycle induced from a proper subgroup. On the other hand, if $G$ is abelian, the usual definition of nondegeneracy is that the skew-symmetric form $\hat G\times\hat G\to\T$, $(\chi,\omega)\mapsto\F(\chi,\omega)\overline{\F(\omega,\chi)}$, is nondegenerate. It is well-known that these are the same conditions. For non-abelian groups there is a similar characterization of nondegeneracy, namely, the following properties are equivalent:
\begin{list}{}{}
\item{(i)} the cocycle $\F$ is nondegenerate;
\item{(ii)} the space $\{(\omega\otimes\iota)(\F_{21}\F^*)\mid \omega\in W^*(G)_*\}$ is weakly operator dense in $W^*(G)$.
\end{list}
This is proved in \cite{Wa1} and in a different way in \cite{La}. The implication (ii)$\Rightarrow$(i) is straightforward, since if $\F=(u\otimes u)(\Ind^G_H\E)\Dhat(u)^*$, then the space in (ii) is contained in $uW^*(H)u^*$. On the other hand, the implication (i)$\Rightarrow$(ii) is quite nontrivial. Below we will sketch a proof of it using ideas of Wassermann \cite[Theorem~12]{Wa1} and of Etingof and Gelaki \cite[Theorem~2.2]{EG0}, \cite[Theorem~3.1]{EG1}.

Let $\RR=\F_{21}\F^*$ and $\Dhat_\F=\F\Dhat(\cdot)\F^*$. Then $(\widehat{\C[G]},\Dhat_\F,\RR)$ is a triangular discrete Hopf  $*$-algebra in the terminology of Appendix~\ref{sb}, where $\widehat{\C[G]}\subset W^*(G)$ is the algebra of matrix coefficients of finite dimensional representations of $G$ with convolution product. By Theorem~\ref{tRadford} and Remark~\ref{rcom}(i) there exists a discrete Hopf $*$-subalgebra $B\subset W^*(G)$ (with respect to the comultiplication $\Dhat_\F$) such that for the von Neumann algebra $N$ generated by $B$, we have $\RR\in N\bar\otimes N$ and the space $\{(\omega\otimes\iota)(\RR)\mid \omega\in W^*(G)_*\}$ is weakly operator dense in $N$.

Consider now the braided C$^*$-tensor category $\CC$ of nondegenerate finite dimensional representations of $B$ with braiding defined by the action of $\Sigma\RR$, where $\Sigma$ is the flip. Since $\RR_{21}\RR=1$, the category $\CC$ is symmetric. It is also even in the sense of \cite[Definition B.8]{Mu}. Indeed, we have a symmetric unitary tensor functor $F\colon \Rep G\to \CC$ which on the level of objects and morphisms is simply the functor of restriction of scalars from~$W^*(G)$ to $B$, while the tensor structure is given by the action of $\F^*$. Since $\Rep G$ is even and every object in $\CC$ is a subobject of an object in the image of $F$, we conclude that~$\CC$ is also even. By the abstract reconstruction theorem of Doplicher and Roberts \cite{DoRo1,DoRo2} (see also \cite{Mu} for an alternative proof), it follows that $\CC$ is equivalent to the symmetric tensor category $\Rep H$ for a compact group~$H$. Composing such an equivalence with the forgetful functor $\CC\to\Hilb$, we get a unitary fiber functor $\tilde F\colon\Rep H\to\Hilb$. This functor is automatically dimension preserving by Corollary~\ref{cdimpres}, hence, ignoring the tensor structure, is isomorphic to the forgetful functor. Fixing such an isomorphism, the tensor structure on $\tilde F$ is defined by $\E^*$ for a cocycle $\E\in Z^2(\hat H;\T)$. It follows that $(N,\Dhat_\F,\F_{21}\F^*)$ is isomorphic to $(W^*(H),\Dhat_\E,\E_{21}\E^*)$, see e.g. \cite[Proposition~1.1]{NT3} for more details. Fix such an isomorphism $\alpha\colon W^*(H)\to N$.

The element $(\alpha\otimes\alpha)(\E^*)\F$ is a symmetric unitary $2$-cocycle on $\hat G$, hence it is the coboundary of a unitary element $u\in W^*(G)$. Consider the cocycle $\G=(u^*\otimes u^*)\F\Dhat(u)$ cohomologous to $\F$ and the $*$-homomorphism $\beta=u^*\alpha(\cdot)u\colon W^*(H)\to W^*(G)$. Since $\F\Dhat(u)=(\alpha\otimes\alpha)(\E)(u\otimes u)$, we have $(\beta\otimes\beta)(\E)=\G$. Since we also have $(\beta\otimes\beta)\Dhat_\E=\Dhat_\G\beta$, it follows that $\beta\colon W^*(H)\to W^*(G)$ respects the untwisted comultiplications, so $\beta$ defines an embedding of $H$ into $G$.

To summarize, we have shown that for any cocycle $\F\in Z^2(\hat G;\T)$ there exists a closed subgroup $H$ of $G$ and a unitary $u\in W^*(G)$ such that $\G=(u^*\otimes u^*)\F\Dhat(u)$  lies in $W^*(H)\bar\otimes W^*(H)$ and the space $\{(\omega\otimes\iota)(\F_{21}\F^*)\mid \omega\in W^*(G)_*\}$ is weakly operator dense in $uW^*(H)u^*$. Assuming now that~$\F$ is nondegenerate, that is, $\F$ is not cohomologous to a cocycle induced from a proper subgroup, we conclude that $H=G$, so the space $\{(\omega\otimes\iota)(\F_{21}\F^*)\mid \omega\in W^*(G)_*\}$ is weakly operator dense in~$W^*(G)$.

\smallskip

Although we will not need this, note that the first part of the above argument shows that any triangular discrete Hopf $*$-algebra
$(A,\Delta,\RR)$ such that its C$^*$-tensor category of nondegenerate finite dimensional representations is even, is isomorphic to $(\widehat{\C[H]},\Dhat_\E,\E_{21}\E^*)$ for a compact group $H$ and a cocycle $\E\in Z^2(\hat H;\T)$. Using the full strength of the Doplicher-Roberts reconstruction theorem we can similarly conclude that any triangular discrete Hopf $*$-algebra $(A,\Delta,\RR)$ is isomorphic to $(\widehat{\C[H]},\Dhat_\E,\RR_k\E_{21}\E^*)$, where $k$ is a central element of order $2$ in $H$ and $\RR_k=(1\otimes1+k\otimes 1+1\otimes k-k\otimes k)/2$. It follows
that the square of the antipode $S$ on $A$ is the identity, see Remark~\ref{rnormalized} above. A direct computation shows then that
the C$^*$-tensor category of nondegenerate finite dimensional representations of $A$ is even if and only if $m(\RR)=1$, where $m\colon M(A\otimes A)\to M(A)$ is the multiplication map;
note that as a byproduct we conclude that $m(\E_{21}\E^*)=1$ for any $\E\in Z^2(\hat H;\T)$, hence $m(\RR_k\E_{21}\E^*)=k$ if $k\in H$ is central of order $2$. The element $m(\RR)$ is nothing else than the Drinfeld element $m(S\otimes\iota)(\RR_{21})$, see \cite[Proposition 4.2.3]{CP}. Indeed, using $S^2=\iota$, $\RR_{21}=\RR^{-1}$ and the identity $(S\otimes\iota)(\RR)=\RR^{-1}$, which holds in any quasitriangular Hopf algebra, we get $\RR=(S\otimes\iota)(\RR_{21})$.

\bigskip

\section{Compact connected groups} \label{s2}

In this section we will apply the general theory presented in the previous section to analyze a particular class of dual cocycles on compact connected groups.

\smallskip

For a compact group $G$ we say that a cocycle $\F\in Z^2(\hat G;\T)$ is invariant if it commutes with elements of the form $\lambda_g\otimes\lambda_g$. Denote by $Z^2_G(\hat G;\T)$ the set of invariant unitary $2$-cocycles. It is easy to see that if $\F$ and $\G$ are cocycles and $\F$ is invariant then $\F\G$ and $\G\F$ are again cocycles. Furthermore, if $\F$ is an invariant cocycle then $\F^*$ is also an invariant cocycle. Therefore $Z^2_G(\hat G;\T)$ is a group.

We say that two invariant cocycles $\F$ and $\G$ are cohomologous if there exists a central unitary element $u\in W^*(G)$ such that $\F=(u\otimes u)\G\Dhat(u)^*$. Denote by $H^2_G(\hat G;\T)$ the group of cohomology classes of invariant unitary $2$-cocycles.

\smallskip

We have an obvious map $H^2_G(\hat G;\T)\to H^2(\hat G;\T)$. To analyze when it is injective, consider the group $\Aut_c(G)$ of automorphisms of $G$ which preserve the inner conjugacy classes. Equivalently, an automorphism $\alpha$ of $G$ belongs to $\Aut_c(G)$ if and only if for any irreducible representation $U$ of $G$ the representations  $U\circ\alpha$ and $U$ are equivalent. The following result is essentially contained in the proof of \cite[Theorem~11]{Wa1}, see also \cite[Corollary~1.8]{GK}.

\begin{proposition} \label{pequivinj}
For a compact group $G$ the following conditions are equivalent:
\enu{i} the map $H^2_G(\hat G;\T)\to H^2(\hat G;\T)$ is injective;
\enu{ii} any symmetric invariant unitary $2$-cocycle on $\hat G$ is the coboundary of  a central unitary element;
\enu{iii} $\Aut_c(G)=\Inn(G)$.
\end{proposition}

\bp (i)$\Leftrightarrow$(ii). Since any symmetric cocycle is a coboundary (of a not necessarily central element), the implication (i)$\Rightarrow$(ii) is immediate. Conversely, assume (ii) holds and $\E$ and $\F$ are invariant cocycles such that $\E=(u\otimes u)\F\Dhat(u)^*$ for a unitary $u$. Then by invariance of $\F$ we have $\E\F^*=(u\otimes u)\Dhat(u)^*$.
Hence, by assumption, $\E\F^*$ is the coboundary of a central unitary element $v$, so $\E$ and $\F$ define the same class in $H^2_G(\hat G;\T)$.

\smallskip

(ii)$\Leftrightarrow$(iii). It is clear that elements of $\Aut_c(G)$ are exactly those automorphisms $\alpha$ of $G$ that can be implemented by unitaries $u\in W^*(G)$, that is, $\lambda_{\alpha(g)}=u\lambda_g u^*$. Furthermore, if $u\in W^*(G)$ is a unitary such that $u^*\lambda(G)u\subset\lambda(G)$, then $u^*\lambda(G)u=\lambda(G)$ by \cite[Lemma~28(2)]{Wa1}. Observe next that $u^*\lambda(G)u\subset\lambda(G)$ if and only if the cocycle $(u\otimes u)\Dhat(u)^*$ is invariant. Indeed, it is invariant if and only if
$$
(\lambda_gu\otimes\lambda_gu)\Dhat(u)^*=(u\otimes u)\Dhat(u)^*(\lambda_g\otimes\lambda_g)
=(u\otimes u)\Dhat(u^*\lambda_g),
$$
equivalently, $u^*\lambda_g u$ is group-like, that is, $u^*\lambda(G)u\subset\lambda(G)$.

Assume now that (ii) holds and $\alpha\in\Aut_c(G)$. Then $\alpha$ is implemented by a unitary $u\in W^*(G)$. The cocycle $(u\otimes u)\Dhat(u)^*$ is invariant, hence $(u\otimes u)\Dhat(u)^*=(v\otimes v)\Dhat(v)^*$ for a central unitary element $v$ by assumption. Then $v^*u$ is group-like, so that $u=v\lambda_g$ for some $g\in G$, and therefore $\alpha=\Ad g$.

Conversely, assume (iii) holds and $\F$ is a symmetric invariant cocycle. Then $\F=(u\otimes u)\Dhat(u)^*$ for a unitary $u$. Since $\F$ is invariant, we have $u^*\lambda(G)u=\lambda(G)$, so $\Ad u$ defines an element $\alpha$ of $\Aut_c(G)$. By assumption, $\alpha=\Ad g$ for some $g\in G$. Then $v=u\lambda_g^*$ is a central unitary and $\F=(v\otimes v)\Dhat(v)^*$.
\ep

The following result will play a key role in the computation of $H^2_G(\hat G;\T)$.

\begin{theorem} \label{tinj}
If $G$ is a compact connected group then the equivalent conditions (i)-(iii) in Proposition \ref{pequivinj} are satisfied.
\end{theorem}

\bp The fact that condition (iii) holds is presumably known to experts. It is e.g.~proved in \cite[Corollary~2]{McMul} and \cite[Lemma~29(2)]{Wa1}, the key point being that any automorphism of a compact connected semisimple Lie group which is trivial on a maximal torus, is inner.
The theorem is also proved by a different method in our paper~\cite{NT4}, where condition (ii) is checked. This method works also for $q$-deformations of compact connected Lie groups. To be more precise, the result in~\cite{NT4} is formulated for simply connected semisimple Lie groups, but it is easy to extend it to arbitrary compact connected groups. Indeed, with minor modifications the proof in~\cite{NT4} works for any compact connected Lie group, basically one needs to replace the lattice of integral weights by the weight lattice. For a general compact connected group $G$ choose a decreasing net of closed normal subgroups $N_i$ such that $\cap_i N_i=\{e\}$ and every quotient $G/N_i$ has a faithful finite dimensional representation, so it is a Lie group. Given a unitary symmetric invariant cocycle $\F$ on $\hat G$ consider the image $\F_i$ of~$\F$ under the homomorphism $W^*(G\times G)\to W^*(G/N_i\times G/N_i)$. Then $\F_i$ is a unitary symmetric invariant cocycle on $\widehat{G/N_i}$, hence the coboundary of a central unitary element $v_i\in W^*(G/N_i)$. Lift $v_i$ to a central unitary element $u_i\in W^*(G)$. Let $u$ be a weak operator limit point of the net $\{u_i\}_i$. Then $u$ is a central unitary element and $\F$ is the coboundary of $u$.
\ep

For finite groups the following result is due to Etingof and Gelaki~\cite[Section~3]{EG2} and Izumi and Kosaki \cite[Lemma~4.1]{IK}. A slightly weaker result was proved by Landstad~\cite[Proposition~4.9]{La}, see also~\cite[3B6]{La0}.

\begin{proposition} \label{pabelian}
Let $G$ be a compact group and $\F\in Z^2(\hat G;\T)$ a cocycle such that the comultiplication $\Dhat_\F=\F\Dhat(\cdot)\F^*$ on $W^*(G)$ is cocommutative. Then there exists a closed normal abelian subgroup $A$ of $G$ and a cocycle $\E\in Z^2(\hat A;\T)$ such that
\enu{i} $\F$ is cohomologous to $\Ind^G_A\E$;
\enu{ii} $\E$ is nondegenerate;
\enu{iii} the class $[\E]$ of $\E$ in $H^2(\hat A;\T)$ is $G$-invariant.

\smallskip\noindent
Furthermore, the pair $(A,[\E])$ is uniquely determined by these conditions. Conversely, if $A$ is a closed normal abelian subgroup of $G$ and $\E\in Z^2(\hat A;\T)$ is such that $[\E]$ is $G$-invariant, then for the cocycle $\F=\Ind^G_A\E$ the comultiplication $\Dhat_\F$ is cocommutative.
\end{proposition}

\bp By Theorem~\ref{tred} there exists a closed subgroup $A$ of $G$ and a cocycle $\E\in Z^2(\hat A;\T)$ such that conditions (i) and (ii) are satisfied. Let $\G=\Ind^G_A\E$. Then the comultiplications $\Dhat_\G\colon W^*(G)\to W^*(G)\otimes W^*(G)$ and  $\Dhat_\E\colon W^*(A)\to W^*(A)\otimes W^*(A)$ are cocommutative.

Consider the $R$-matrix $\RR=\E_{21}\E^*$ for the comultiplication $\Dhat_\E$, so in particular we have the identity $(\Dhat_\E\otimes\iota)(\RR)=\RR_{13}\RR_{23}$. By nondegeneracy of $\E$ the space $\{(\omega\otimes\iota)(\RR)\mid \omega\in W^*(A)_*\}$ is weakly operator dense in $ W^*(A)$. Using cocommutativity of $\Dhat_\E$, we have
$$
\RR_{23}\RR_{13}=(\sigma\otimes\iota)(\RR_{13}\RR_{23})=(\sigma\otimes\iota)(\Dhat_\E\otimes\iota)(\RR)
=(\Dhat_\E\otimes\iota)(\RR)=\RR_{13}\RR_{23},
$$
where $\sigma$ is the flip. It follows that $\{(\omega\otimes\iota)(\RR)\mid \omega\in W^*(A)_*\}$ is a commutative algebra; see also Remark~\ref{rcom}(ii) below. Hence the group $A$ is abelian.

By cocommutativity of $\Dhat_\G$ the element $\E^*\E_{21}$ commutes with $\Dhat(\lambda_g)$ for every $g\in G$. Since~$W^*(A)$ is commutative, we have $\E^*\E_{21}=\E_{21}\E^*=\RR$, so $(\Ad\lambda_g\otimes\Ad\lambda_g)(\RR)=\RR$. Using again nondegeneracy of $\E$ we see that the subalgebra $W^*(A)\subset W^*(G)$ is $\Ad\lambda_g$-invariant for every $g\in G$. Hence $A$ is normal in $G$.

Since $A$ is abelian, two cocycles $\mathcal A$ and $\mathcal B$ in $Z^2(\hat A;\T)$ are cohomologous if and only if $\mathcal A\mathcal B^*$ is symmetric, that is, $\mathcal A_{21}\mathcal A^*=\mathcal B_{21}\mathcal B^*$. Since $\RR=\E_{21}\E^*$ is $G$-invariant, we conclude that the cocycles~$\E$ and $(\Ad\lambda_g\otimes\Ad\lambda_g)(\E)$ are cohomologous, that is, $[\E]\in H^2(\hat A;\T)$ is $G$-invariant.

Finally, the uniqueness statement follows from that in Theorem~\ref{tred}.

\smallskip

Conversely, if $\F=\Ind^G_A\E$, $A$ is normal abelian and $[\E]$ is $G$-invariant, then $\E_{21}\E^*$ commutes with $\lambda_g\otimes\lambda_g$ by $G$-invariance of $[\E]$. Since $\E_{21}\E^*=\E^*\E_{21}$, it follows that $\Dhat_\F$ is cocommutative.
\ep

Since for invariant cocycles $\F$ we have $\Dhat_\F=\Dhat$, it follows that we get a map from $H^2_G(\hat G;\T)$ into the set of pairs $(A,[\E])$, where $A$ is a closed normal abelian subgroup of $G$ and $[\E]\in H^2(\hat A;\T)^G$. For finite groups this map has been analyzed by Guillot and Kassel in \cite[Section~4]{GK}. Rather than trying to generalize their analysis to compact groups, we will be content with the following theorem, which is the principal result of this section.

\begin{theorem} \label{tmain0}
For any compact connected group $G$ we have
$$
H^2_G(\hat G;\T)\cong H^2(\widehat{Z(G)};\T),
$$
where $Z(G)$ is the center of $G$.
\end{theorem}

\bp Assume $\F\in Z^2_G(\hat G;\T)$. Then the twisted coproduct $\Dhat_\F$ coincides with $\Dhat$. Hence by Proposition~\ref{pabelian} the cocycle $\F$ is cohomologous to a cocycle $\G$ induced from a closed normal abelian subgroup $A\subset G$. Since $G$ is connected, $A$ is contained in the center of $G$, because the action of $G$ on the discrete group $\hat A$ can only be trivial. Thus $\G$ is induced from a cocycle on $\widehat{Z(G)}$. But any such cocycle
is invariant. Hence, by Proposition~\ref{pequivinj} and Theorem~\ref{tinj}, $\F$ and $\G$ are cohomologous as invariant cocycles. Therefore the homomorphism $H^2(\widehat{Z(G)};\T)\to H^2_G(\hat G;\T)$, $\E\mapsto\Ind^G_{Z(G)}\E$, is surjective. It is clearly injective, since a cocycle on $\widehat{Z(G)}$ is a coboundary if it is symmetric, and any coboundary on $\hat G$ is symmetric.
\ep

For a C$^*$-tensor category $\CC$ denote by $\Aut^\otimes(\CC)$ the group of unitary $\C$-linear monoidal autoequivalences of $\CC$ identified up to unitary monoidal natural isomorphism. Denote by $\Out(G)$ the group $\Aut(G)/\Inn(G)$ of outer automorphisms of $G$.

\begin{theorem} \label{tmain}
For any compact connected group $G$ we have a canonical isomorphism
$$
\Aut^\otimes(\Rep G)\cong H^2(\widehat{Z(G)};\T)\rtimes \Out(G).
$$
\end{theorem}

\bp We define a homomorphism $\gamma\colon H^2(\widehat{Z(G)};\T)\rtimes \Out(G)\to \Aut^\otimes(\Rep G)$ as follows.
For $\alpha\in\Aut(G)$ define an autoequivalence $\tilde\alpha$ of $\Rep G$ by letting $\tilde\alpha(\pi_U)=\pi_U\circ\alpha^{-1}$ for every unitary representation $\pi_U\colon G\to B(H_U)$, while the action of $\tilde\alpha$ on morphisms, as well as the tensor structure of $\tilde\alpha$, are defined in the obvious way. Then denoting by $[\alpha]$ the class of $\alpha$ in $\Out(G)$, put~$\gamma([\alpha])$ to be the class of $\tilde\alpha$. On the other hand, for a unitary cocycle $\E$ on $\widehat{Z(G)}$, define an autoequivalence $\beta_\E$ which acts trivially on objects and morphisms, while the tensor structure is given by~$\E^*$ considered as an element of $W^*(G)\bar\otimes W^*(G)$. Then put $\gamma([\E])$ to be the class of $\beta_\E$. It is easy to see that $\gamma$ is well-defined and $\tilde\alpha\circ\beta_\E=\beta_{(\alpha\otimes\alpha)(\E)}\circ\tilde\alpha$, so $\gamma$ is indeed a homomorphism.

If $\tilde\alpha\circ\beta_\E$ is naturally isomorphic to the identity functor then $\alpha$ acts trivially on the isomorphism classes of irreducible representations, so $\alpha\in\Aut_c(G)$. By Theorem~\ref{tinj} we have $\alpha\in\Inn(G)$. It follows that $\beta_\E$ is naturally isomorphic to the identity functor, which is the same as saying that $\Ind^G_{Z(G)}\E$ is a coboundary of a central element in $W^*(G)$. But the map $H^2(\widehat{Z(G)};\T)\to H^2_G(\hat G;\T)$ is injective, so $\E$ is a coboundary. Thus $\gamma$ is injective.

To prove surjectivity, consider a monoidal autoequivalence $F$ of $\Rep G$. It defines an automorphism of the fusion ring of $G$ mapping irreducibles into irreducibles. By \cite[Theorem~2]{McMul} any such automorphism is implemented by an automorphism of $G$. Thus replacing $F$ by $\tilde\alpha\circ F$ for an appropriate $\alpha$ we may assume that $F$ maps every representation to an equivalent one. Since $\Rep G$ is semisimple, $F$ is monoidally isomorphic to a functor which acts trivially on objects and morphisms. The  tensor structure of the latter functor is defined by an invariant cocycle $\F$ on $\hat G$. By Theorem~\ref{tmain0} the cocycle $\F$ is cohomologous (as an invariant cocycle) to an induced cocycle $\Ind^G_{Z(G)}\E$. Hence $F$ is monoidally isomorphic to $\beta_\E$.
\ep

If $\CC$ is a symmetric C$^*$-tensor category, denote by $\Aut^\otimes_s(\CC)\subset\Aut^\otimes(\CC)$ the subgroup of symmetric autoequivalences of $\CC$. Note that in the notation of the above proof, the functor $\tilde\alpha\circ\beta_\E$ is symmetric if and only if $\E$ is symmetric, hence a coboundary of a central element. Therefore $\Aut^\otimes_s(\Rep G)\cong \Out(G)$ for a compact connected group $G$. This result, however, follows from the general fact that a compact group is completely determined by the symmetric tensor category of its finite dimensional representations, see e.g.~\cite[Theorem B.6]{Mu}, which implies that $\Aut^\otimes_s(\Rep G)\cong \Aut(G)/\Aut_c(G)$ for any compact group $G$. In particular, if $G$ is a simply connected semisimple compact Lie group then $\Aut^\otimes_s(\Rep G)$ is isomorphic  to the group of automorphisms of the based root datum of $G$. 

\smallskip

Although a compact group is determined by its symmetric C$^*$-tensor category of finite dimensional representations, this is not so if one forgets the symmetric structure. Two compact groups~$G_1$ and~$G_2$  are called monoidally equivalent if the C$^*$-tensor categories $\Rep G_1$ and $\Rep G_2$ are equivalent. In the next section we will extend a result of Etingof and Gelaki~\cite{EG2} and Izumi and Kosaki~\cite{IK} and describe all pairs of monoidally equivalent compact separable groups. The situation is much simpler when one of the groups is connected.

\begin{theorem} \label{tmonod}
Assume $G$ is a compact connected group. Then any compact group which is mono\-idally equivalent to $G$, is isomorphic to $G$.
\end{theorem}

\bp Assume a compact group $\tilde G$ is monoidally equivalent to $G$, so there exists a monoidal equivalence $E\colon \Rep G\to\Rep\tilde G$. Consider the forgetful functors $F\colon\Rep G\to\Hilb$ and $\tilde F\colon\Rep\tilde G\to\Hilb$. Since $E$ is intrinsic dimension preserving, the fiber functor $\tilde F\circ E\colon \Rep G\to\Hilb$ is dimension preserving, hence, ignoring the tensor structure, it is naturally isomorphic to $F$. Fixing such a unitary isomorphism, the tensor structure of $\tilde F\circ E$ is defined by $\F^*$, where $\F$ is a unitary $2$-cocycle on $\hat G$. Identifying $W^*(\tilde G)$ with the subalgebra of natural transformations of $\tilde F\circ E$, we conclude that the Hopf-von Neumann algebra of $\tilde G$ is isomorphic to $(W^*(G),\Dhat_\F)$. The comultiplication $\Dhat_\F$ is cocommutative, hence by Proposition~\ref{pabelian}, $\F$ is cohomologous to a cocycle $\G$ induced from a closed normal abelian subgroup of $G$. But since $G$ is connected, such a subgroup is contained in the center of $G$, so that $\Dhat_\G=\Dhat$. We thus see that the Hopf-von Neumann algebra of $\tilde G$ is isomorphic to $(W^*(G),\Dhat)$, hence $\tilde G\cong G$.
\ep

In fact, the result of McMullen~\cite{McMul} implies a stronger result: if the fusion rings of $G$ and $\tilde G$ are isomorphic, via an isomorphism which maps irreducibles into irreducibles, and $G$ is connected, then~$G$ and $\tilde G$ are isomorphic. This is the case since isomorphism of fusion rings implies that $\tilde G$ is connected as well, which, in turn, is true because a compact group is connected if and only if tensor powers of any nontrivial irreducible representation contain infinitely many mutually inequivalent irreducible subrepresentations, see e.g.~\cite[Theorem~28.21]{HR}. The point of Theorem~\ref{tmonod}, however, is that as opposed to McMullen's result it can be proved without using any structure theory of compact groups.

\bigskip

\section{Monoidally equivalent groups} \label{s3}

In this section we will generalize results of Etingof and Gelaki \cite[Theorem~1.3]{EG2} and Izumi and Kosaki~\cite[Section~4]{IK} and describe all pairs of monoidally equivalent compact separable groups.

\smallskip

Let $K$ be a compact group acting by automorphisms on a compact abelian group $A$; let $g.a$ denote the action of $g\in K$ on $a\in A$. Assume that both groups $K$ and $A$ are separable. Following~\cite{EG2} we will define a homomorphism
$$
\tau\colon H^2(\hat A;\T)^K\to H^2(K;A),
$$
where $H^2(K;A)$ denotes Moore cohomology with Borel cochains~\cite{Moore3}. Let $\tilde c\in Z^2(\hat A;\T)$ be a cocycle such that it cohomology class $c=[\tilde c]$ is $K$-invariant. For $g\in K$ denote by $\tilde c^g$ the cocycle defined by $\tilde c^g(\varphi,\chi)=\tilde c(g^{-1}\varphi,g^{-1}\chi)$, where $(g^{-1}\varphi)(a)=\varphi(g.a)$ for $a\in A$. By assumption the cocycles~$\tilde c$ and~$\tilde c^g$ are cohomologous, so there exists a $1$-cochain $z_g\in C^1(\hat A;\T)$ such that $\tilde c=\tilde c^gdz_g$. Furthermore, the group $C^1(\hat A;\T)$ is compact and we may assume that the map $K\ni g\mapsto z_g\in C^1(\hat A;\T)$ is Borel. Indeed,
the subset
$$
C:=\{(g,z)\in K\times C^1(\hat A;\T)\mid \tilde c=\tilde c^gdz\}\subset K\times C^1(\hat A;\T)
$$
is closed, so the projection $C\to K$ onto the first coordinate has a Borel section~\cite{Kal}. For $g,h\in K$ define
$$
\tilde b(g,h)=z_g(z_h)^gz_{gh}^{-1}\in C^1(\hat A;\T),
$$
where $(z_h)^g(\varphi)=z_h(g^{-1}\varphi)$. Then $\tilde b$ is a $C^1(\hat A;\T)$-valued Borel $2$-cocycle on $K$. We have
$$
d(\tilde b(g,h))=dz_g(dz_h)^g(dz_{gh})^{-1}=\tilde c(\tilde c^g)^{-1}\tilde c^g(\tilde c^{gh})^{-1}\tilde c^{-1}\tilde c^{gh}=1,
$$
so that $\tilde b(g,h)\in A\subset C^1(\hat A;\T)$. It is easy to check that the class $b=[\tilde b]$ of $\tilde b$ in $H^2(K;A)$ depends only on the class $c\in H^2(\hat A;\T)^K$. We put $\tau(c)=b$.

\smallskip

It is clear that if the cocycle $\tilde c$ is itself $K$-invariant, then $\tilde b$ is trivial. In other words, the kernel of~$\tau$ contains the image of the canonical homomorphism $H^2_K(\hat A;\T)\to H^2(\hat A;\T)$.

\smallskip

For $b\in H^2(K;A)$ denote by $G_b$ the corresponding extension of $K$ by $A$. Recall~\cite{Mac1,Mac2,Moore3} that one first defines $G_b$ as a Borel group and then shows that it caries a unique topology making it a compact group.

\begin{theorem}
Assume $K$ and $A$ are compact separable groups, $A$ is abelian, and $K$ acts on $A$ by automorphisms. Assume $b_1,b_2\in H^2(K;A)$ are such that $b_2b_1^{-1}\in \tau(H^2(\hat A;\T)^K)$. Then the groups~$G_{b_1}$ and~$G_{b_2}$ are monoidally equivalent. Any pair of monoidally equivalent compact separable groups is obtained this way for appropriate $K$ and $A$.
\end{theorem}

\bp Let $G_1$ and $G_2$ be monoidally equivalent compact separable groups. As we have already observed in the proof of Theorem~\ref{tmonod},  we can then identify $(W^*(G_2),\Dhat_2)$ with $(W^*(G_1),(\Dhat_1)_\F)$ for some $\F\in Z^2(\hat G_1;\T)$. By Proposition~\ref{pabelian} we may further assume that $\F=\Ind^{G_1}_A\E$ for a closed normal abelian subgroup $A\subset G_1$ and a nondegenerate cocycle $\E\in Z^2(\hat A;\T)$ such that the cohomology class of $\E$ is $G_1$-invariant.

Put $K=G_1/A$. Let $b_1\in H^2(K;A)$ be the element corresponding to the extension $G_1$ of $K$ by~$A$. Choose a Borel section $s_1\colon K\to G_1$ of the quotient map with $s_1(e)=e$, and let $\tilde b_1(g,h)=s_1(g)s_1( h)s_1(gh)^{-1}$, so $\tilde b_1$ is a cocycle representing $b_1$.

Since $[\E]\in H^2(\hat A;\T)^{G_1}=H^2(\hat A;\T)^K$, we can apply the above procedure and construct a cocycle $\tilde b\in Z^2(K;A)$ such that $\tilde b(g, h)=z_{ g}(z_{h})^{g}z_{gh}^{-1}$ and $\E=\E^{g}dz_{g}$; we may assume that $z_e=1$. We can think of $z_g\in C^1(\hat A;\T)$ as a unitary element of $W^*(A)\subset W^*(G_1)$. To simplify the notation we will also identify $G_1$ with the group of group-like elements of $(W^*(G_1),\Dhat_1)$. Then we can write
\begin{equation}\label{eiso}
{\tilde b(g, h)}z_{gh}=z_{ g}{s_1(g)}z_{h}{s_1(g)}^{-1}\ \ \hbox{in}\ \ W^*(G_1)
\end{equation}
and
$$
\E=({s_1(g)}\otimes {s_1(g)})\E({s_1(g)}\otimes {s_1(g)})^{-1}(z_{g}\otimes z_g)\Dhat_1(z_g)^{-1}
=(z_g{s_1(g)}\otimes z_g{s_1(g)})\E\Dhat_1(z_g{s_1(g)})^{-1}.
$$
The last identity means that $z_gs_1(g)$ is a group-like element for the twisted coproduct $(\Dhat_1)_\E=\Dhat_2$, so it is an element of $G_2\subset W^*(G_1)$. We therefore get a Borel map $s_2\colon K\to G_2$, $s_2(g)=z_gs_1(g)$. Note also that since $\Dhat_2=\Dhat_1$ on $W^*(A)$, the group $A$ is a subgroup of $G_2$. Consider the cocycle $\tilde b_2=\tilde b\tilde b_1$ and put $b_2=[\tilde b_2]\in H^2(K;A)$. Define the group $G_{b_2}$ using the cocycle $\tilde b_2$. We have a Borel map $\varphi\colon G_{b_2}\to G_2$, $\varphi(ag)=as_2(g)$ for $a\in A$ and $g\in K$.
By virtue of~\eqref{eiso} it is a group homomorphism. We claim that it is an isomorphism, hence an isomorphism of topological groups (see e.g.~\cite[Theorem~9.10]{Kec}).

Let us check first that $\varphi$ is injective. Assume $\varphi(ag)=1\in W^*(G_1)$, that is, $az_gs_1(g)=1$. Then $z_g\in G_1\cap W^*(A)=A$. Hence $g=e$, $z_g=z_e=1$ and $a=e$.

To prove surjectivity note that the above arguments imply that for every $u\in G_1$ there exists $z\in W^*(A)$ such that $zu\in G_2$. But the roles of $G_1$ and $G_2$ are completely symmetric: since $(W^*(G_1),\Dhat_1)=(W^*(G_2),(\Dhat_2)_{\F^*})$, by the proof of Proposition~\ref{pabelian} we conclude that $A$ is normal in~$G_2$ and the cohomology class $[\E]$ is $G_2$-invariant. Hence for every $v\in G_2$ there exists $z\in W^*(A)$ such that $zv\in G_1$. Assume $zv=as_1(g)$. Then $zz_gv=az_gs_1(g)=\varphi(ag)$. It follows that $zz_g\in G_2\cap W^*(A)=A$, and therefore $v=\varphi((zz_g)^{-1}ag)$.

Thus we have proved that $G_1\cong G_{b_1}$, $G_2\cong G_{b_2}$ and $b_2b_1^{-1}=\tau([\E])$.

\smallskip

Conversely, assume $G_1\cong G_{b_1}$, $G_2\cong G_{b_2}$ and $b_2b_1^{-1}=\tau([\E])$ for some $\E\in Z^2(\hat A;\T)$ such that~$[\E]$ is $K$-invariant. Note that we do not assume that $\E$ is nondegenerate, but we may assume that $\E=\Ind^A_B\E_0$ for a uniquely defined $K$-invariant subgroup $B\subset A$ and a nondegenerate cocycle $\E_0\in Z^2(\hat B;\T)$, namely, $B\subset A$ is the annihilator of the kernel of the skew-symmetric form $\E_{21}\E^*$. Fix a Borel section $s_1\colon K\to G_1$ and let $\tilde b_1\in Z^2(K;A)$ be the corresponding cocycle, so $[\tilde b_1]=b_1$. Put $\F=\Ind^{G_1}_A\E$. Consider the Hopf-von Neumann algebra $(W^*(G_1),(\Dhat_1)_\F)$. By Proposition~\ref{pabelian} invariance of $[\E]$ implies that $(\Dhat_1)_\F$ is cocommutative. By the Tannaka-Krein reconstruction we conclude that $(W^*(G_1),(\Dhat_1)_\F)=(W^*(G),\Dhat)$ for a compact group $G$. The groups $G_1$ and $G$ are monoidally equivalent, so to complete the proof it suffices to show that $G_2\cong G$.

As above, let $\tilde b$ be a cocycle representing $\tau([\E])$, so $\tilde b(g, h)=z_{ g}(z_{h})^{g}z_{gh}^{-1}$. Consider the cocycle $\tilde b_2=\tilde b\tilde b_1$. Then $\tilde b_2$ corresponds to a Borel section $s_2\colon K\to G_2$. We have a Borel homomorphism $\varphi\colon G_2\to G\subset W^*(G_1)$, $\varphi(as_2(g))=az_gs_1(g)$. The same argument as before shows that $\varphi$ is injective, while to prove surjectivity it suffices to show that for every $v\in G$ there exists $z\in W^*(A)$ such that $zv\in G_1$. Furthermore, we know that the latter property would be satisfied if $\E$ were nondegenerate. But by assumption $\E=\Ind^A_B\E_0$ and $\E_0$ is nondegenerate, so for every $v\in G$ there exists $z\in W^*(B)\subset W^*(A)$ such that $zv\in G_1$. Hence $\varphi$ is an isomorphism.
\ep

Examples of nonisomorphic monoidally equivalent finite groups can be found in \cite{EG2} and \cite{IK}.

\bigskip

\appendix

\section{Amenability} \label{sappa}

Let $\CC$ be a $C^*$-tensor category. Throughout the whole section we will assume that $\CC$ is small, strict and has subobjects, conjugates and irreducible unit $\unit$, consult e.g.~\cite{LR} for definitions. 

A dimension function on $\CC$ is a map
$d \colon {\rm Ob} \, {\mathcal C} \to [1,+\infty)$ such that
$$
d(U) = d(V)\ \ \hbox{if}\ \ U \cong V,\ \ d(U \oplus V) = d(U) + d(V),\ \  d(U\otimes V) = d(U)d(V),\ \
d(\bar U) = d(U).
$$
Consider the fusion ring $K({\mathcal C})$ of $\CC$. Fix representatives $U_j$, $j \in I$, of isomorphism classes of simple objects. Denote by $e\in I$ the point corresponding to $\unit$. Then we can think of $K(\CC)$ as $\oplus_{j \in I}{\mathbb Z}j$, with multiplication given by
$$
ij = \sum_k m_{ij}^kk,\ \ \hbox{if}\ \ U_i \otimes U_j \simeq \bigoplus_k m_{ij}^k \, U_k.
$$
Define an involution $j\mapsto\bar j$ on $I$ by  $U_{\bar j}\cong \bar U_j$. Then we can say that a dimension function is a homomorphism $d \colon K({\mathcal C}) \to {\mathbb R}$ such that
$$
d(j) \geq 1 ,\ \ d(\bar j) = d(j).
$$

There always exists at least one dimension function, namely, the intrinsic dimension $d_i$, see~\cite{LR}. Let us briefy recall how it is defined. An object $\bar U$ is called a conjugate of an object $U$ if there exist morphisms $R\colon\unit\to \bar U\otimes U$ and $\bar R\colon\unit\to U\otimes\bar U$ such that the compositions
$$
U\xrightarrow{\iota\otimes R} U\otimes\bar U\otimes U\xrightarrow{\bar R^*\otimes\iota} U\ \ \hbox{and}\ \
\bar U\xrightarrow{\iota\otimes \bar R} \bar U\otimes U\otimes\bar U\xrightarrow{R^*\otimes\iota}\bar U
$$
are the identity maps. By assumption every object has a conjugate object. The latter is unique up to isomorphism. The intrinsic dimension of $U$ is then defined as
$$
d_i(U)=\min\{\|R\|\cdot\|\bar R\|\},
$$
where the minimum is taken over all possible morphisms $R$ and $\bar R$ as above.

\smallskip

For every $j\in I$ denote by $\Lambda_j$ the operator of multiplication by $j$ on $\C\otimes_\Z K(\CC)=\oplus_{k\in I}\C k$. The following lemma, applied to $\gamma_{ik}=m^i_{jk}/d(j)=m^k_{\bar j i}/d(j)$ and $s_i=t_i=d(i)$, shows that it extends to a bounded linear operator on $\ell^2(I)$ and $\|\Lambda_j\|\le d(j)$, where $d$ is any dimension function on $\CC$.

\begin{lemma} \label{lcontraction}
Let $\Gamma=(\gamma_{ik})_{i,k\in I}$ be a matrix with nonnegative real coefficients. Assume there exists a vector $s=(s_i)_{i\in I}$ such that $s_i>0$ for all $i\in I$, the vector $t=\Gamma s$ is well-defined and $\Gamma't\le s$ coordinate-wise, where $\Gamma'$ is the transposed matrix. Then $\Gamma$ defines a contraction on~$\ell^2(I)$.
\end{lemma}

\bp Note that if $t_i=0$ for some $i$ then $\gamma_{ik}=0$ for all $k$. Thus replacing zero by any strictly positive number for every such $i$ we get a vector $t$ such that $t_i>0$ for all $i$, $\Gamma s\le t$ and $\Gamma't\le s$. Then for vectors $\xi=(\xi_i)_i$ and $\zeta=(\zeta_i)_i$ in $\ell^2(I)$, using the Cauchy-Schwarz inequality we get
\begin{align*}
|(\Gamma\xi,\zeta)|&=\left|\sum_{i,k}\gamma_{ik}\xi_k\bar\zeta_i\right|
=\left|\sum_{i,k}\left((\gamma_{ik}t_is_k^{-1})^{1/2}\xi_k\right)
\left((\gamma_{ik}s_kt_i^{-1})^{1/2}\bar\zeta_i\right)\right|\\
&\le\left(\sum_{i,k}\gamma_{ik}t_is_k^{-1}|\xi_k|^2\right)^{1/2}
\left(\sum_{i,k}\gamma_{ik}s_kt_i^{-1}|\zeta_i|^2\right)^{1/2}\le\|\xi\|\,\|\zeta\|,
\end{align*}
so that $\|\Gamma\|\le1$.
\ep

Fix now a dimension function $d$ on $\CC$. For a probability measure $\mu$ on $I$ define a contraction $\lambda_\mu$ on~$\ell^2(I)$ by
$$
\lambda_\mu=\sum_{j\in I}\frac{\mu(j)}{d(j)}\Lambda_j.
$$
We will write $\lambda_j$ instead of $\lambda_{\delta_j}$. If $\mu$ and $\nu$ are two probability measures then $\lambda_{\mu}\lambda_{\nu}=\lambda_{\mu*\nu}$, where
$$
(\mu*\nu)(k)=\sum_{i,j}m_{ i  j}^k\frac{d(k)}{d(i)d(j)} \mu(i)\nu(j).
$$
We will write $\mu^n$ for the $n$-th convolution power $\mu*\dots*\mu$ of $\mu$. For a measure $\mu$ denote by $\check{\mu}$ the measure defined by $\check{\mu}(i)=\mu(\bar i)$. A probability measure $\mu$ on $I$ is called symmetric if $\check{\mu}=\mu$, and it is called nondegenerate if $\cup_{n\ge1}\operatorname{supp}\mu^{n}=I$.

\smallskip

The following lemma is contained in \cite[Theorem~4.1 and Proposition~4.8]{HI} and is an analogue of known results for random walks on groups.

\begin{lemma}
For a probability measure $\mu$ on $I$ consider the following conditions:
\enu{i} $1\in\operatorname{Sp}\lambda_\mu$;
\enu{ii} $\|\lambda_\mu\|=1$;
\enu{iii}  $(\check{\mu}*\mu)^{n}(e)^{1/n}\to1$ as $n\to+\infty$.

\smallskip\noindent
Then $(\mathrm i)\Rightarrow(\mathrm {ii})\Leftrightarrow(\mathrm {iii})$. If $\mu$ is symmetric then all three conditions are equivalent, and if they are satisfied, there exists a sequence $\{\xi_n\}_n$ of unit vectors in $\ell^2(I)$ such that $\|\lambda_{j}\xi_n-\xi_n\|\to0$ as $n\to+\infty$ for all $j\in\cup_{n\ge1}\operatorname{supp}\mu^n$; in particular, $1\in\Sp\lambda_{\nu}$ for any probability measure $\nu$ such that $\operatorname{supp}\nu$ is contained in $\cup_{n\ge1}\operatorname{supp}\mu^n$.
\end{lemma}

\bp It is clear that (i)$\Rightarrow$(ii), and since $(\check{\mu}*\mu)^{n}(e)=((\lambda_\mu^*\lambda_\mu)^n\delta_e,\delta_e)$, that (iii)$\Rightarrow$(ii).

To show that (ii)$\Rightarrow$(iii), consider the unital C$^*$-algebra $A$ generated by $\lambda_\mu^*\lambda_\mu$. By Lemma~\ref{lcontraction} the operators of multiplication by $j\in I$ on the right on $K(\CC)$ extend to bounded operators on $\ell^2(I)$. It follows that the vector $\delta_e$ is cyclic for the commutant $A'$ of $A$, hence $(\cdot\,\delta_e,\delta_e)$ is a faithful state on~$A$. From this we conclude that $\|\lambda_\mu\|^2$ is the least upper bound of the support of the measure $\nu$ on $\Sp\lambda_\mu^*\lambda_\mu$ defined by the state $(\cdot\,\delta_e,\delta_e)$. Since
$$
((\lambda_\mu^*\lambda_\mu)^n\delta_e,\delta_e)=\int t^n\,d\nu(t)\ \ \hbox{for all}\ \ n\ge0,
$$
it is easy to see that this upper bound is equal to
$$
\lim_{n\to+\infty}((\lambda_\mu^*\lambda_\mu)^n\delta_e,\delta_e)^{1/n}
=\lim_{n\to+\infty}(\check{\mu}*\mu)^{n}(e)^{1/n}.
$$
Hence (ii)$\Rightarrow$(iii).

\smallskip

Next assume that $\mu$ is symmetric and condition (ii) is satisfied. Since $\lambda_\mu$ is self-adjoint and $\|\lambda_\mu\|=1$, there exists a sequence of unit vectors $\zeta_n\in\ell^2(I)$ such that $|(\lambda_\mu\zeta_n,\zeta_n)|\to1$. Consider the unit vectors $\xi_n=|\zeta_n|$. Since the matrix $\lambda_\mu$ has nonnegative coefficients, we have $(\lambda_\mu\xi_n,\xi_n)\ge|(\lambda_\mu\zeta_n,\zeta_n)|$. Hence $(\lambda_\mu\xi_n,\xi_n)\to1$, and therefore $1\in\Sp\lambda_\mu$.

Since $\lambda_\mu$ is a convex combination of the operators $\lambda_j$, we also see that $(\lambda_j\xi_n,\xi_n)\to1$ for every $j\in\operatorname{supp}\mu$. Since $\lambda_j$ is a contraction, this implies that $\|\lambda_{j}\xi_n-\xi_n\|\to0$. Since $\|\lambda^k_\mu\xi_n-\xi_n\|\to0$ for all $k\ge1$, we similarly conclude that $\|\lambda_{j}\xi_n-\xi_n\|\to0$ for every $j\in\operatorname{supp}\mu^k$.
\ep

We thus see that for a fixed dimension function $d$ on $\CC$ the following conditions are equivalent:
\begin{list}{}{}
\item{(i)} $1\in\operatorname{Sp}\lambda_\mu$ for every probability measure $\mu$;
\item{(ii)} $\|\lambda_\mu\|=1$ for every probability measure $\mu$;
\item{(iii)}  $(\check{\mu}*\mu)^{n}(e)^{1/n}\to1$ as $n\to+\infty$ for every probability measure $\mu$.
\end{list}
Furthermore, it suffices to check any of these conditions for a net $\{\mu_\alpha\}_\alpha$ of symmetric probability measures such that the sets $\cup_{n\ge1}\operatorname{supp}\mu^{n}_\alpha$ increase with $\alpha$ and their union is $I$; in particular, if $I$ is countable, it suffices to consider one nondegenerate symmetric probability measure. If the above conditions are satisfied, following Hiai and Izumi~\cite{HI} we say that the pair $(K(\CC),d)$ is amenable. See \cite[Section~4]{HI} for other equivalent characterizations of amenability.

\smallskip

We will say that the category $\CC$ is amenable if $(K(\CC),d_i)$ is amenable, where $d_i$ is the intrinsic dimension function.
For an object $U\cong\oplus_jn_jU_j$ consider the probability measure $\mu=d_i(U)^{-1}\sum_jn_jd_i(U_j)\delta_j$. The condition $\|\lambda_\mu\|=1$ was discussed by Longo and Roberts in~\cite[Section~3]{LR} as a possible definition of amenability of the object $U$. Thus $\CC$ is amenable if and only if every object is amenable.

\begin{proposition}
Let $d$ be a dimension function on a C$^*$-tensor category $\CC$ such that $(K(\CC),d)$ is amenable. Then $d(j)=\|\Lambda_j\|=\|(m^i_{jk})_{i,k}\|$ for every $j\in I$, and for any other dimension function $d'$ on $\CC$ we have $d'\ge d$.
\end{proposition}

\bp
For $(K(\CC),d)$ to be amenable, at the very least we need the condition $\|\lambda_{j}\|=1$ to be satisfied for every $j\in I$, but it means exactly that $d(j)=\|\Lambda_j\|$. The second statement follows from the inequality $\|\Lambda_j\|\le d'(j)$, which holds for any dimension function $d'$.
\ep

In view of the above proposition it is natural to ask whether the map $j\mapsto\|\Lambda_j\|$ always defines a dimension function; see \cite[Theorem~8.2]{ENO} for a solution of a similar problem for fusion categories.

\smallskip

If $\CC$ is finite (that is, $I$ is finite) then any dimension function is amenable, since the vector $(d(j))_{j\in I}\in\ell^2(I)$ is an eigenvector of $\lambda_\mu$ with eigenvalue $1$ for every probability measure $\mu$. It follows that $\CC$ has only one dimension function, the intrinsic dimension. This can also be easily seen from the Perron-Frobenius theorem~\cite{LR,ENO}. Without the finiteness assumption there is at least the following result, mentioned in \cite[Section~3]{LR}.

\begin{proposition} \label{pdimpres0}
Let $\CC$ and $\CC'$ be C$^*$-tensor categories and $F\colon\CC\to\CC'$ be a unitary tensor functor. Assume $\CC$ is amenable. Then $d_i(F(U))=d_i(U)$ for every object $U$ in $\CC$.
\end{proposition}

\bp Consider the dimension function $d=d_i\circ F$ on $\CC$. By definition of intrinsic dimension it is clear that $d\le d_i$. On the other hand, $d_i$ is the smallest dimension function on $\CC$ by amenability. Hence $d=d_i$.
\ep

We will next explore the relation between amenability of fusion rings and amenability of quantum groups. Let $G$ be a compact quantum group~\cite{Wor}, and $(\C[G],\Delta)$ be the Hopf $*$-algebra of matrix coefficients of finite dimensional corepresentations of $G$. 
Denote by $C_u(G)$ the universal enveloping C$^*$-algebra of $\C[G]$. Let $h$ be the Haar state on $C_u(G)$, and $\pi_h\colon C_u(G)\to B(L^2(G))$ be the corresponding GNS-representation with cyclic vector $\xi_h$. Put $C_r(G)=\pi_h(C_u(G))$.

\smallskip

The following result is essentially due to Banica \cite[Section~6]{Ba1}, with the key idea attributed in \cite{Ba1} to Skandalis; see also~\cite{BMT}.

\begin{theorem} \label{tamen}
For a compact quantum group $G$ the following conditions are equivalent:
\enu{i} the fusion ring of $G$ with the classical dimension function $\dim$ is amenable;
\enu{ii} the counit $\eps\colon\C[G] \to {\mathbb C}$ extends to a character of $C_r (G)$;
\enu{iii} the map $\pi_h\colon C_u (G) \to C_r (G)$ is an isomorphism.
\end{theorem}

If these conditions are satisfied, the quantum group $G$ is called coamenable. See \cite{Tom}, \cite{BCT} and references therein for other equivalent definitions. We stress that this is not the same as amenability of the category $\Rep G$ of finite dimensional corepresentations of $G$, since coamenability refers to classical dimension and not intrinsic dimension. The two notions of dimension coincide if and only if $G$ is of Kac type.

\begin{proof}[Proof of Theorem~\ref{tamen}] We will first give an equivalent formulation of condition (i). Let $U\in B(H_U)\otimes\C[G]$ be a finite dimensional unitary corepresentation of $G$: $(\iota\otimes\Delta)(U)=U_{12}U_{13}$. Define the character of $U$ by $\chi (U) = (\Tr \otimes \pi_h) (U)\in C_r(G)$. We claim that (i) is equivalent to $\dim U\in\Sp\chi(U)$ for every~$U$.

To see this, consider the closure $A$ of the linear span of the operators $\chi(U)$ in $C_r(G)$. This is a unital C$^*$-algebra: $\chi(U)^*=\chi(\bar U)$ and $\chi(U)\chi(V)=\chi(U\times V)$, where $U\times V=U_{13}V_{23}$ is the tensor product of the corepresentations $U$ and $V$. Put $H=\overline{A\xi_h}\subset L^2(G)$. Since the Haar state is faithful on $C_r(G)$, we have $\dim U\in\Sp\chi(U)$ if and only if $\dim U\in\Sp(\chi(U)|_H)$.

Choose representatives $U_j$, $j\in I$, of isomorphism classes of irreducible unitary corepresentations of~$G$. By the orthogonality relations the map $T\colon\ell^2(I)\to H$, $T\delta_j=\chi(U_j)\xi_h$, is a unitary isomorphism. Furthermore, if  $U\cong\oplus_jn_jU_j$  and $\mu=(\dim U)^{-1}\sum_jn_j(\dim U_j)\delta_j$ then $T^*\chi(U)T=(\dim U)\lambda_\mu$. Therefore $1\in\Sp\lambda_\mu$ if and only if $\dim U\in\Sp\chi(U)$, and our equivalent formulation of~(i) is proved.

\smallskip

(i)$\Leftrightarrow$(ii). Assume (ii) holds. Then $\eps(\chi(U))=\dim U$ as $(\iota\otimes\eps)(U)=1$, and since $\eps$ is a character of the C$^*$-algebra $C_r(G)$, it follows that $\dim U\in \Sp\chi(U)$.

Assume (i) holds. Let $U$ be a self-conjugate unitary corepresentation. Since $\dim U\in\Sp\chi(U)$ and~$\chi(U)$ is self-adjoint, there exists a state $\omega$ on $C_r(G)$ such that $\omega(\chi(U))=\dim U$. Since the inequality $\|\chi(V)\|\le\dim V$ holds for any corepresentation $V$, we conclude that $\omega(\chi(V))=\dim V$ for any subcorepresentation $V$ of $U$. Choosing an increasing net of self-conjugate unitary corepresentations and taking a weak$^*$ limit point of the corresponding states~$\omega$, we get a state $\nu$ on~$C_r(G)$ such that $\nu(\chi(V))=\dim V$ for any $V$. For every~$V$, the operator $X=(\iota\otimes\nu)(V)\in B(H_V)$ has the properties $\|X\|\le1$ and $\Tr X=\dim H_V$, which is possible only when $X=1$. Hence $\nu$ is an extension of the counit $\eps$ on $\C[G]$.

\smallskip

(ii)$\Leftrightarrow$(iii). The implication (iii)$\Rightarrow$(ii) is obvious, as the counit is well-defined on $C_u(G)$. To prove that (ii)$\Rightarrow$(iii), observe that the comultiplication $\Delta\colon\C[G]\to \C[G]\otimes\C[G]$ extends to a $*$-homomorphism $\alpha\colon C_r(G)\to C_r(G)\otimes C_u(G)$. Indeed, $\alpha$ is implemented by the right regular corepresentation $W\in M( K(L^2 (G) \otimes C_u(G)))$ of $G$, that is, $\alpha(a)=W(a\otimes1)W^*$ for $a\in C_r(G)$ and~$W$ is defined by $W(\pi_h(b)\xi_h\otimes c)=(\pi_h\otimes\iota)(\Delta(b))(\xi_h\otimes c)$ for $b\in \C[G]$ and $c\in C_u(G)$. Then $\rho=(\eps\otimes\iota)\alpha\colon C_r(G)\to C_u(G)$ is a unital $*$-homomorphism such that $\rho(\pi_h(a))=a$ for $a\in \C[G]$. Hence $\rho$ is the inverse of $\pi_h\colon C_u(G)\to C_r(G)$.
\ep

\begin{corollary} \label{cdimpres}
Any compact group $G$ is coamenable. In particular, if $F\colon \Rep G\to\Hilb$ is a unitary fiber functor, then $\dim F(U)=\dim U$ for any $U$.
\end{corollary}

\bp Clearly, $C_r(G)=C(G)$ and the counit $\eps\colon C(G)\to\C$ is well-defined: it is the evaluation at $e\in G$. Therefore $\Rep G$ is amenable, so from Proposition~\ref{pdimpres0} we get the second statement.
\ep

The statement that any unitary fiber functor on $\Rep G$ is dimension preserving can also be deduced by other methods. For example, from the ergodic actions point of view it is equivalent to the fact that the notions of full multiplicity and full quantum multiplicity~\cite{BDRV} coincide for compact groups. Without the unitarity assumption the result is not true already for $G=\SU(2)$~\cite{Bru}. Note also that by \cite[Example~7.6]{HI} for compact Lie groups one has a stronger result than amenability of~$\Rep G$: the fusion ring of $G$ has polynomial growth.

\smallskip

Turning to genuine quantum groups we get the following result~\cite{Ba1,Ba2}.

\begin{corollary}
The Drinfeld-Jimbo $q$-deformation $G_q$ ($q>0$) of any simply connected semisimple compact Lie group~$G$ is coamenable.
\end{corollary}

\bp For $q=1$ the result is true by the previous corollary. On the other hand, it is well-known that the fusion ring of $G_q$ together with the classical dimension function do not depend on $q$.
\ep

\bigskip

\section{Minimal Hopf subalgebras} \label{sb}

The goal of this appendix is to prove an analogue of Radford's theorem \cite{Rad} on minimal Hopf subalgebras of quasitriangular Hopf algebras for compact quantum groups. It will be more convenient to formulate the result in the equivalent language of discrete Hopf $*$-algebras \cite{VD,NT3}. Recall briefly that this means that we consider $*$-algebras $A$ of the form $\oplus_jB(H_j)$ and the comultiplication is a nondegenerate $*$-homomorphism $\Delta\colon A\to M(A\otimes A)=\prod_{j,k}B(H_j)\otimes B(H_k)$. A discrete Hopf $*$-algebra $A$ is called quasitriangular if we are given an invertible element $\RR\in M(A\otimes A)$ such that
$$
\Delta^{op} =\RR\Delta(\cdot )\RR^{-1},\ \
(\Delta\otimes\iota)(\RR)=\RR_{13}\RR_{23},\ \
(\iota\otimes\Delta)(\RR)=\RR_{13}\RR_{12}\ \ \hbox{and}\ \ \RR^*=\RR_{21}.
$$
Denote by $\hat A$ the dual Hopf $*$-algebra $\oplus_jB(H_j)^*\subset M(A)^*$ with product $(\omega\eta)(a)=(\omega\otimes\eta)\Delta(a)$, comultiplication $\Dhat(\omega)(a\otimes b)=\omega(ab)$ and involution $\omega^*(a)=\overline{\omega(S(a)^*)}$, where $S$ is the antipode on~$A$.

\begin{theorem} \label{tRadford}
Let $A$ be a quasitriangular discrete Hopf $*$-algebra with $R$-matrix $\RR$. Then there exists a discrete Hopf $*$-subalgebra $B\subset M(A)$ sitting in $M(A)$ nondegenerately (that is, $BA=A$) such that $\RR\in M(B\otimes B)$ and the linear span $\U_\RR$ of elements of the form $(\omega\otimes\iota)(\RR)(\iota\otimes\eta)(\RR)$, $\omega,\eta\in \hat A$, is a $\sigma(M(A),\hat A)$-dense Hopf $*$-subalgebra of $M(B)$.
\end{theorem}

\bp We start by proving that $\U_\RR$ is a Hopf $*$-algebra. The identities $(\Delta\otimes\iota)(\RR)=\RR_{13}\RR_{23}$, $(\iota\otimes\Delta)(\RR)=\RR_{13}\RR_{12}$ and $(S\otimes S)(\RR)=\RR$ show that the map $\hat A^{cop}\ni\omega\mapsto(\omega\otimes\iota)(\RR)\in M(A)$ is a Hopf algebra homomorphism, so its image $\U_l$ is a Hopf algebra. Similarly the map $\hat A^{op}\ni\eta\mapsto(\iota\otimes\eta)(\RR)\in M(A)$ is a Hopf algebra homomorphism, so its image $\U_r$ is also a Hopf algebra. Therefore to show that $\U_\RR=\U_l\U_r$ is a Hopf algebra we just have to check that $\U_r\U_l\subset \U_l\U_r$. This is a consequence of the Yang-Baxter equation $\RR_{12}\RR_{13}\RR_{23}=\RR_{23}\RR_{13}\RR_{12}$. Indeed, denote by $\tilde\RR$ the $R$-matrix $\RR$ considered as an element of $M(A\otimes A^{op})$. Then the Yang-Baxter equation can be written as
$$
\RR_{12}\tilde\RR_{23}\tilde\RR_{13}=\tilde\RR_{13}\tilde\RR_{23}\RR_{12}\ \ \hbox{in}\ \ M(A\otimes A\otimes A^{op}),
$$
and flipping the second and the third factor we get
$$
\RR_{13}\tilde\RR_{32}\tilde\RR_{12}=\tilde\RR_{12}\tilde\RR_{32}\RR_{13}\ \ \hbox{in}\ \ M(A\otimes A^{op}\otimes A).
$$
Observe next that since $(\iota\otimes S^{-1})(\RR)=\RR^{-1}$ and $S^{-1} \colon A^{op}\to A$ is an algebra isomorphism, the element $\tilde \RR$ is invertible with inverse $(\iota\otimes S)(\tilde\RR)$. Therefore we can write
$\tilde\RR_{32}\RR_{13}=\tilde\RR_{12}^{-1}\RR_{13}\tilde\RR_{32}\tilde\RR_{12}$. Applying $\omega\otimes\eta\otimes\iota$ to this identity and writing $(\omega\otimes\eta)(\tilde\RR_{12}^{-1}\cdot\tilde\RR_{12})$ as $\sum_i\omega_i\otimes\eta_i$, we get
\begin{align*}
(\iota\otimes\eta)(\RR)(\omega\otimes\iota)(\RR)
&=(\omega\otimes\eta\otimes\iota)(\tilde\RR_{32}\RR_{13})
=\sum_i(\omega_i\otimes\eta_i\otimes\iota)(\RR_{13}\tilde\RR_{32})\\
&=\sum_i(\omega_i\otimes\iota)(\RR)(\iota\otimes\eta_i)(\RR).
\end{align*}
Thus $\U_\RR$ is a Hopf algebra. Finally, to see that it is a $*$-algebra, note that since $\RR^*=\RR_{21}$, we have $x\in\U_l$ if and only if $x^*\in\U_r$.

\smallskip

Let $M\subset M(A)$ be the von Neumann algebra generated by $A$, that is, if we identify the $*$-algebra~$A$ with $\oplus_{j\in I}B(H_j)$, then $M$ is the $\ell^\infty$-direct sum of the C$^*$-algebras $B(H_j)$. Let $N\subset M$ be the von Neumann subalgebra consisting of elements $x\in M$ such that $\pi(x)\in\pi(\U_\RR)$ for every nondegenerate finite dimensional representation $\pi\colon A\to B(K)$. Then $\pi(N)=\pi(\U_\RR)$ for every $\pi$. Indeed,
for every finite subset $\alpha\subset I$ consider the representation $\pi_\alpha$ of $A$ on $\oplus_{j\in\alpha}H_j$. Let $\beta$ be such that $\pi$ is quasi-equivalent to $\pi_\beta$. For every $\alpha\supset\beta$ denote by $\rho_\alpha$ the projection $\oplus_{j\in\alpha}B(H_j)\to\oplus_{j\in\beta}B(H_j)$. Let $a\in\pi_\beta(\U_\RR)$. For every $\alpha\supset\beta$ choose $a_\alpha\in \pi_\alpha(\U_\RR)$ such that $\rho_\alpha(a_\alpha)=a$ and $\|a_\alpha\|=\|a\|$. Then choose $x_\alpha\in M$ such that $\pi_\alpha(x_\alpha)=a_\alpha$ and $\|x_\alpha\|=\|a_\alpha\|=\|a\|$. Take a weak operator limit point $x$ of the net $\{x_\alpha\}_\alpha$. Then $x\in N$ and $\pi_\beta(x)=a$.

Since $M$ is a finite discrete von Neumann algebra, $N$ is also finite and discrete, so $N$ is the $\ell^\infty$-direct sum of full matrix algebras $B_k$. Put $B=\oplus_kB_k\subset N$. We claim that $\Delta(B)\subset N\bar\otimes N$. To show this, it suffices to check that $(\pi_\alpha\otimes\pi_\beta)\Delta(N)\subset \pi_\alpha(\U_\RR)\otimes \pi_\beta(\U_\RR)$ for any finite $\alpha,\beta\subset I$. The representation $(\pi_\alpha\otimes\pi_\beta)\Delta$ is quasi-equivalent to $\pi_\gamma$ for some $\gamma$. Since $\pi_\gamma(N)=\pi_\gamma(\U_\RR)$, it follows that
$(\pi_\alpha\otimes\pi_\beta)\Delta(N)=(\pi_\alpha\otimes\pi_\beta)\Delta(\U_\RR)$, and since $\Delta(\U_\RR)\subset \U_\RR\otimes\U_\RR$, our claim is proved. Similarly one checks that $S(B)=B$. Hence $B$ is indeed a discrete Hopf $*$-subalgebra of $M(A)$. It sits in~$M(A)$ nondegenerately, since $N$ contains the unit of $M$.

Since $\pi_\alpha(N)=\pi_\alpha(\U_\RR)$ for all finite $\alpha$, it is clear that $\U_\RR$ is contained in the algebra of closed densely defined operators affiliated with $N$, that is, $\U_\RR\subset M(B)$, and $\U_\RR$ is $\sigma(M(A),\hat A)$-dense in~$M(B)$. Finally, since $(\pi_\alpha\otimes\pi_\beta)(\RR)\in \pi_\alpha(\U_\RR)\otimes\pi_\beta(\U_\RR)$ for all finite $\alpha$ and $\beta$, the operator $\RR$ is affiliated with~$N\bar\otimes N$, so $\RR\in M(B\otimes B)$.
\ep

\begin{remark} \label{rcom} \mbox{\ }
\enu{i} If $\RR$ is unitary, in which case we say that $A$ is a triangular discrete Hopf $*$-algebra, then the identities $(S\otimes\iota)(\RR)=\RR^*=\RR_{21}$ imply that $\U_l=\U_r$, so in this case
$$
\U_\RR=\{(\omega\otimes\iota)(\RR)\mid \omega\in\hat A\}=\{(\iota\otimes\eta)(\RR)\mid\eta\in \hat A\},
$$
and the homomorphisms $\hat A^{cop}\ni\omega\mapsto(\omega\otimes\iota)(\RR)\in M(A)$ and $\hat A^{op}\ni\eta\mapsto(\iota\otimes\eta)(\RR)\in M(A)$ are $*$-preserving.
\enu{ii} If $A$ is commutative or cocommutative then $B$ is both commutative and cocommutative, so $B\cong\widehat{\C[G]}$ for a compact abelian group $G$. Indeed, if $A$ is commutative (resp., cocommutative) then $\hat A$ is cocommutative (resp., commutative), so that $\U_l$ and $\U_r$ are both commutative and cocommutative; see also the proof of Proposition~\ref{pabelian}. Therefore we just have to show that if $A$ is cocommutative then~$\U_l$ and $\U_r$ commute with each other. To see this, note that commutativity of $\U_l$ and $\U_r$ means that $\RR_{13}$ commutes with $\RR_{12}$ and $\RR_{23}$. The Yang-Baxter equation $\RR_{12}\RR_{13}\RR_{23}=\RR_{23}\RR_{13}\RR_{12}$ implies then that $\RR_{12}$ and $\RR_{23}$ also commute. This means exactly that $\U_l$ and $\U_r$ commute with each other.
\end{remark}

\end{document}